\documentclass[12pt]{amsart}

\usepackage{todonotes}

\usepackage{amsmath}
\usepackage{amsthm}
\usepackage{amssymb}
\usepackage[colorlinks]{hyperref}

\usepackage[nodisplayskipstretch]{setspace}

\makeatletter
\newcommand{\reqnomode}{\tagsleft@false}
\makeatother

\renewcommand{\subset}{\subseteq}

\newcommand{\kk}{k}
\newcommand{\polynomial}{\kk[t]}
\newcommand{\rational}{\kk(t)}
\newcommand{\cusp}{\kk[t^2,t^3]}
\newcommand{\laurent}{\kk[t,t^{-1}]}

\newcommand{\oa}{\mathsf{a}}
\newcommand{\ob}{\mathsf{b}}
\newcommand{\oc}{\mathsf{c}}
\newcommand{\oK}{\mathsf{K}}
\newcommand{\oC}{\mathsf{C}}
\newcommand{\oD}{\mathsf{D}}
\newcommand{\oE}{\mathsf{E}}
\newcommand{\oF}{\mathsf{F}}
\newcommand{\oL}{\mathsf{L}}

\newcommand{\oR}{\mathsf{R}}
\newcommand{\oX}{\mathsf{X}}
\newcommand{\oY}{\mathsf{Y}}
\newcommand{\oZ}{\mathsf{Z}}

\newcommand{\thetab}{\theta}
\newcommand{\thetah}{\theta}

\newcommand{\GL}{\mathrm{GL}}
\newcommand{\GLn}{\mathrm{GL}_n}

\newcommand{\Mtwo}{\mathrm{M}_2}
\newcommand{\Mthree}{\mathrm{M}_3}
\newcommand{\Msix}{\mathrm{M}_6}
\newcommand{\Mts}{\mathrm{M}_{36}}
\newcommand{\Mn}{\mathrm{M}_n}
\newcommand{\Mnm}{\mathrm{M}_{n-1}}

\newcommand{\cop}{\mathrm{cop}}
\newcommand{\op}{\mathrm{op}}
\newcommand{\ev}{\mathrm{v}}

\newcommand{\act}{\triangleright}

\newcommand{\ls}{[\![}
\newcommand{\rs}{]\!]}

\newtheorem*{theoremm}{Theorem}

\newtheorem{lemma}{Lemma}[subsection]
\newtheorem{proposition}{Proposition}[subsection]
\newtheorem{corollary}{Corollary}[subsection]

\theoremstyle{definition}

\newtheorem{definition}{Definition}[subsection]

\theoremstyle{remark}

\newtheorem{remark}{Remark}[subsection]
\newtheorem{example}{Example}[subsection]

\begin{document}
\title{Pointed Hopf algebra
(co)actions on
rational functions}
\author{Ulrich Krähmer}
\author{Blessing Bisola Oni}
\email{ulrich.kraehmer@tu-dresden.de}
\email{blessing@aims.edu.gh}

\begin{abstract}
This article studies 
the construction of Hopf algebras $H$
acting on a given algebra $K$ in terms 
of algebra morphisms $ \sigma \colon K \rightarrow
\Mn(K)$. The approach is
particularly suited for controlling whether
these actions restrict to a given
subalgebra $B$ of $K$, whether $H$ is
pointed, and whether these actions are
compatible with a given $*$-structure on
$K$. The theory is applied to the field 
$K=\rational$ of rational functions 
containing the coordinate ring 
$B=\cusp$ of the cusp. An explicit example
is described in detail
and shown to define a new quantum
homogeneous space structure on
the cusp.
\end{abstract}
\maketitle

\tableofcontents

\section{Introduction}
The question which Hopf algebras can (co)act
on a commutative algebra goes back
at least to Cohen's article
\cite{cohen}.
For semisimple Hopf
algebras $H$ over algebraically closed fields
$k$ of characteristic $0$, it was
answered completely by Etingof and Walton 
\cite{etingofwalton}, see also
\cite{eitngofwaltongoswami}: if $H$ acts
inner faithfully on a commutative integral domain
$K$ (i.e.~the action does not arise from
the action of a proper quotient Hopf
algebra of $H$), then $H$ is a 
group algebra. In contrast to
this, there are typically many inner
faithful actions of pointed Hopf
algebras on a commutative
algebra $K$, even 
when $K$ admits only few automorphisms and
derivations. 
For example, the coordinate rings $B$ of singular plane
curves are conjecturally all 
quantum homogeneous spaces in the sense of 
\cite{mullerschneider}, that is, can be
embedded as right coideal subalgebras into
Hopf algebras $A$ which are faithfully
flat as modules over these subalgebras
\cite{ulan,ulma,anke}, and in
many examples studied so far, these
embeddings yield
inner faithful actions of pointed Hopf algebras $H$ that are
dually paired with $A$.

The starting point of the present article
is to demand that like classical
symmetries (given by actions of group algebras or
universal enveloping algebras), the
action of $H$ extends
from the coordinate ring $B$ to the field
$K$ of rational functions on the curve. 
As the classification of Hopf algebra
(co)actions on fields is an
interesting topic in its own right (see
e.g.~\cite{etingofwaltonpointed1,etingofwaltonpointed2, 
greitherpareigis, taylortruman, Tsang,cresporiovela}
and the references therein), we felt it worthwhile
to begin a systematic study of such
actions on function fields that restrict
to coordinate rings. 

The approach we take
was maybe first applied by Manin
in his construction of
quantum $SL(2)$ as a Hopf algebra
(co)acting on the quantum plane \cite{manin}. More
recently, it was used mostly in the
$C^*$-algebraic quantum group community in
the construction and study of compact
quantum automorphism groups such as the
quantum permutation groups or the
liberations of compact Lie groups
\cite{banicabichon, banicagoswami, weber, 
banicaspeicher, banicabichoncollins, 
bhowmickskalskisoltan,banicabhowmickdecommer}.

In this approach, a bialgebra action on
a $k$-algebra $K$ is given by -- or more
precisely constructed from -- an algebra
morphism $K \rightarrow \Mn(K)$; the
bialgebra is a Hopf algebra if this
morphism, viewed as an element of 
$ \Mn ( \mathrm{End}_k(K))$ is strongly
invertible 
in the sense of
Definition~\ref{StrInv} below. The general
approach is well-known, but some aspects
of our presentation are to our knowledge
novel, such
as the connection to the theory of general
linear groups over noncommutative rings
(see Section~\ref{strinvmat}), the 
treatment of $*$-structures in
the pointed rather than the semisimple setting
(see Section~\ref{starstructures}), and
the application to the
field $K=\rational$ of rational functions (see
Section~\ref{rationalfunctionsec}).  

Our main focus is the
construction of pointed Hopf algebra
(co)actions on $K=\rational$ which restrict to
the coordinate ring
$B=\cusp$ of the cusp. 
Working over a field $k$ whose
characteristic is not $2$ or $3$, we construct Hopf
algebras $A_\sigma$ and $H_\sigma$ which
are generated as algebras by 
$\oa,\ob,\oc$ 
respectively $K, D, Y$. These generators
satisfy the defining relations  
\begin{gather*}
	\oa \ob + \ob \oa =  
	\oa \oc + \oc \oa = 
	\ob \oc + 
	\oc \ob = 0,\\
	3 \ob^2 =
	\oa^6, \quad \oc^2=1,\\
	K D+ DK=KY-YK=0, \quad
 Y^2 D - 2 Y D Y + D Y^2
=0,\\ 
	K^2=1,\quad D^2=0.
\end{gather*} 
In terms of these generators, 
the coproduct, counit, and antipode of
$A_\sigma $ and $H_\sigma $ are given by
\begin{gather*}
	\Delta(\oa) = 1 \otimes \oa
+ \oa \otimes \oc,\quad
	\Delta (\oc) = \oc \otimes 
	\oc \\
	\Delta(\ob) = 1 \otimes \ob + 
	\oa^2 \otimes \oa - \oa
\otimes \oa^2 \oc + \ob
\otimes \oc \\
	\varepsilon(\oa) = 
	\varepsilon(\ob) = 0, \quad
	\varepsilon (\oc) = 1,\\
	S(\oa) = - \oa \oc, \quad 
	S(\ob) = - \ob \oc,\quad
	S(\oc) = \oc,\\
\Delta(K) = K \otimes K, \quad \Delta(D) = 1 \otimes D + D \otimes K, \\
\Delta(Y) = 1 \otimes Y - 6D \otimes DK + Y \otimes 1,\\
\varepsilon(K) = 1, \quad \varepsilon(D) = \varepsilon(Y) = 0, \\
S(K) = K, \quad S(D) = -DK, \quad S(Y) = -Y.
\end{gather*}

Our main results are summarised in the following theorem:

\begin{theoremm}
We have:
\begin{enumerate}
\item The coalgebras
$A_\sigma$ and $H_\sigma$ are pointed.
\item There is a dense 
Hopf algebra embedding 
$ A_\sigma \rightarrow H_\sigma^\circ$. 
\item For any point $(\lambda ^2,\lambda
^3)$, $ \lambda \in k$, of the cusp, there
is an embedding 
$$ 
	\iota \colon B = \cusp
	\rightarrow A_\sigma
$$
of 
its coordinate ring  
as a right coideal subalgebra such that 
$$
	t^2 \mapsto 
	\lambda^2 1 + 
	12 \oa^2,
	\quad 
	t^3 \mapsto 
	6\lambda^2 \oa +
	36 \oa^3 +  	
	36 \ob +
	\lambda^3 \oc 
$$
and $A_\sigma$ is faithfully flat over
$B$. 
\item The 
$H_\sigma$-action on $B$ dual to the
resulting $A_\sigma$-coaction on $B$ is inner faithful.
\item This $H_\sigma$-action and
this $A_\sigma$-coaction both
extend to the field
$K=\rational$ of rational functions. 
\item If $k= \mathbb{C} $, then
$A_\sigma$ and $H_\sigma $ become Hopf
$*$-algebras with 
$$
	K^*  = K,\quad
	D^* = - D, \quad 
	Y^* = -Y + 6 iD,
$$
and the images of $A_\sigma$ and of $B$ 
in $H_\sigma^\circ$ are
$*$-subalgebras if $ \bar \lambda = \lambda $;
the resulting $*$-structure on $
\mathbb{C} [t^2,t^3]$ is
given by $t^*=t$. 
\end{enumerate}
\end{theoremm}

The paper is divided into two main
sections: our presentation is arranged in
such a way that large parts of the development
of the general
theory and of the computations
carried out for the
cusp only use linear algebra and the
theory of polynomials and rational
functions in one variable. No Hopf
algebra theory or algebraic
geometry is required for this
material, which is gathered in 
Section~\ref{theory}. 
Section~\ref{classificationsection}
contains the first steps towards a
classification of pointed Hopf
algebra actions on
$\cusp$ that extend to $\rational$, but
we mostly focus on the specific
example described in the theorem above.  
Section~\ref{interpretation} contains the
proof of this theorem.\\  

\noindent {\bf Acknowledgements.}  
U.K.~is supported
by the DFG grant
``Cocommutative comonoids''
(KR 5036/2-1); 
B.B.O.~is supported
by SISSA Trieste and by ICTP
Trieste.

\section{Quantum automorphisms
of the cusp}\label{theory}
This section
contains those parts of the 
paper that can be 
formulated in terms
of elementary algebra. 
The interpretation
in terms of Hopf algebras
is contained in 
Section~\ref{interpretation}. 

\subsection{Strongly invertible
matrices}\label{strinvmat}

Recall that if 
$
	\sigma \in \GLn(P)
$
is an invertible 
$n \times n$-matrix 
with entries in a (unital associative) 
ring $P$, then the transpose $ \sigma ^T$
is in general not invertible:
\begin{example}
If $a,d \in P$ satisfy 
$da = 1 \neq ad$, then 
we have
$$
	\left(\begin{array}{cc}
	a & 1 \\
	0 & d
\end{array}\right)
	\in \GL_2(P),\quad
	\left(\begin{array}{cc}
	a & 0 \\
	1 & d
	\end{array}\right)
	\notin \GL_2(P).
$$
\end{example}

More
precisely, $ \sigma^T \in \GLn(P)$ if and
only if $ \sigma \in \GLn(P^\op)$, 
where $P^\op$ denotes the
opposite ring of
$P$ (the additive abelian
group $P$ equipped with the opposite multiplication
$a \cdot_\op b:=ba$). Indeed, this follows 
from the straightforwardly verified fact
that for  
$ \sigma , \tau \in \Mn(P)$, we have 
\begin{equation}\label{galoiskunda3}
	(\sigma \tau )^T = \tau ^T \cdot_\op 
	\sigma ^T,
\end{equation}
where on the right hand side, $\cdot_\op$ is the
multiplication in $\Mn(P^\op)$. 
This yields the \emph{contragredient
isomorphism} of general linear groups 
$$
	\GLn(P) \rightarrow 
	\GLn(P^\op),\quad
	\sigma \mapsto 
	\hat \sigma := (\sigma ^{-1})^T,
$$ 
see e.g.~\cite[Chapter~3]{hahnomeara} for a
discussion of this map.

The theory of Hopf algebras motivates to
study those matrices to 
which one can apply
this operation arbitrarily
often without leaving $\GLn(P)$; we are
not aware of a standard name for such
matrices, so we introduce a working
terminology:

\begin{definition}\label{StrInv}
We call $ \sigma \in \Mn(P)$  
\emph{strongly invertible} if
there exists a sequence 
$\{\sigma _d\}_{d \in \mathbb{Z}
}$ in $\GLn(P)$ with
$ \sigma _0 = \sigma $ and 
$\sigma _{d+1} = \hat \sigma _d$. 
\end{definition}

Note this means that also all
$ \sigma _d^T$ are invertible
with  
$(\sigma_d^T)^{-1}= \sigma
_{d-1}$. 
Note also that the strongly
invertible matrices do not
form a group:

\begin{example}\label{exampletoronto}
Let $\kk$ be a field and 
$P:=\kk \langle x,y \rangle $
be the free $\kk$-algebra
generated by $x,y$. Then 
$$
	\sigma :=
	\begin{pmatrix}
	1 & x \\
	0 & 1
	\end{pmatrix},
	\quad
	\tau :=
	\begin{pmatrix}
	1 & 0 \\
	y & 1
	\end{pmatrix}
$$
are strongly invertible
matrices, but 
$ (\sigma \tau)^T $ is not
invertible. 
\end{example}

We furthermore remark that the set of
strongly invertible matrices
is in general a proper subset
of  
$\GL_n(P) \cap \GL_n(P^\op)$:

\begin{example}
Assume that $P$ is any ring that contains
elements $x,y$ and that $x$ is 
invertible.
Consider
$$
	\sigma :=
	\begin{pmatrix}
	x^{-1} + y^2 & y \\
	y & 1
	\end{pmatrix} \in 
	\GL_2(P), \quad 
	\sigma ^{-1} = 
	\begin{pmatrix}
	x & -xy \\
	-yx & 1+yxy
	\end{pmatrix}. 
$$
Since $ \sigma = \sigma ^T$, 
we have $ \sigma \in 
\GL_2(P) \cap \GL_2(P^\op)$. 
However, 
$$ 
	\hat \sigma = (\sigma ^{-1})^T
	= \begin{pmatrix}
	x & -yx \\
	-xy & 1+yxy
	\end{pmatrix}
\in \Mtwo(P)
$$
is in general not invertible: indeed,
assume that 
	$ \hat\sigma \in
\GL_2(P)$ with 
$$
	\begin{pmatrix}
	a & b \\
	c & d 
	\end{pmatrix}
	:= \hat\sigma^{-1}. 
$$
Then the condition  
$ \hat\sigma  
\hat\sigma ^{-1}= 
1$ implies
$$
	b = x^{-1}yxd,\quad
	1 = (1 + y (xy-yx))d. 
$$
Similarly, 
$\hat\sigma^{-1} 
\hat\sigma = 1$
yields  
$$
	c = dxyx^{-1},\quad 
	1 = d
	(1 + yxy -xyx^{-1} yx).
$$
Thus $d \in P$ is invertible 
with inverse
$$
	d^{-1} = 
	1 + yxy -xyx^{-1} yx
	= 1 + y (xy-yx). 
$$
However, this implies 
$(y-xyx^{-1})y=0$ which does
not hold in general.  
\end{example}

In this paper, we will focus on upper
triangular matrices, and for these, the
condition of strong
invertibility 
is easily controlled:

\begin{proposition}\label{StrInv2}
If $ \sigma \in \Mn(P)$ is upper
triangular, $\sigma _{ij} = 0$ for 
$i>j$, then the following
statements are
equivalent:
\begin{enumerate}
\item $ \sigma $ is strongly
invertible.
\item $ \sigma $ is
invertible. 
\item $ \sigma_{ii} \in P$ is invertible for
$i=1,\ldots,n$. 
\end{enumerate}
In this case, 
$ \sigma ^{-1} $ is upper triangular and
all $ (\sigma ^{-1})_{ij}$ are contained
in the subring of $P$ generated by
the $ \sigma_{ij}$ and the $ \sigma _{ii}^{-1}$.  
\end{proposition}
\begin{proof}
``$ \Rightarrow $'':
Suppose $\sigma$ is strongly
invertible. Then 
$ \sigma
$ and $ \sigma ^T$ are
invertible.  
As $ \sigma _{ni}=0$ for $i<n$, we then have 
\begin{equation}\label{toronto}
	(\sigma \sigma^{-1})_{nn}  =
	\sigma_{nn}(\sigma^{-1})_{nn} = 1 
\end{equation}
and similarly
$$
	1 =((\sigma^T)^{-1} \sigma^T)_{nn} =
	((\sigma^T)^{-1})_{nn} \sigma_{nn}. 
$$ 
Hence, $ \sigma _{nn}$ is
invertible in $P$ with inverse
$$
	(\sigma_{nn})^{-1}= (\sigma^{-1})_{nn} = 
	((\sigma^T)^{-1})_{nn}.
$$ 
$\sigma \sigma^{-1}$ shows
$$
	 \sigma_{nn}
	(\sigma^{-1})_{nj} = 0 \;\;
	\forall \; j = 1, \cdots, n-1. 
$$ 
But $\sigma_{nn}$ is invertible, 
thus 
$$
	(\sigma^{-1})_{nj} = 0 \;\; 
	\forall \; j = 1, \cdots,
n-1.
$$ 
Analogously, one shows
that
$ ((\sigma ^T)^{-1})_{jn}$
vanishes 
for $j=1,\ldots, n-1$. 
So $ \sigma $,
$ \sigma ^{-1}$, and 
$ (\sigma ^T)^{-1}$ 
can be written
in block matrix form as 
$$ 
	\sigma = 
	\begin{pmatrix}
\alpha & \mu \\
0^T & \sigma_{nn}
\end{pmatrix}, \qquad
\sigma^{-1} = \begin{pmatrix}
\beta & \gamma \\
0^T & \sigma^{-1}_{nn}
\end{pmatrix},\qquad
	(\sigma ^T)^{-1} = 
	\begin{pmatrix} 
	\delta & 0 \\
	 \nu^T & \sigma _{nn}^{-1}
	\end{pmatrix},
$$
where $\alpha,\beta, \delta \in \Mnm(P)$, $\mu,\gamma,\nu$ are 
column vectors in $P^{n-1}$, 
and $0$ is the zero vector in
$P^{n-1}$.

From $ \sigma \sigma ^{-1} = 1 =
\sigma ^{-1} \sigma $ we
obtain that $ \alpha $ is
invertible with inverse $
\beta $.
Analogously, $ \alpha ^T$ is
invertible with inverse $
\delta $. We also obtain 
$$
	\gamma = - \alpha^{-1} 
	\mu \sigma _{nn}^{-1},
$$
so the entries $ \gamma_j$ are
elements of the subring of $P$
generated by the entries of 
$ \sigma $ and of $
\alpha^{-1} $. 
Continuing inductively one
obtains the claim. 

``$ \Leftarrow $'': 
Suppose the diagonal entries 
$\sigma_{ii}$ of $\sigma = \sigma_0$ are invertible. 
We show that the equation 
$ \sigma \tau = 1$ can be
solved inductively in the ring
of upper triangular matrices
with entries in $P$. Indeed,
this equation means that 
$$
	\sum_{m=l}^n 
	\sigma _{lm} \tau_{mn} =
\delta _{ln}. 
$$
These equations can be solved
by induction on $n-l$. For
$n-l=0$ we obtain 
$$
	\tau _{ll} = \sigma
_{ll}^{-1}.
$$
For $n-l = i$, we obtain $$\tau_{l\; i+l} = - \sigma_{ll}^{-1} 
\bigg( \sum_{m = l +1}^{n} \sigma_{lm}\tau_{mn} \bigg).$$
Furthermore, solving the equation  $\tau \sigma = 1$ 
inductively as above, we conclude that $\tau = \sigma^{-1}$.
Analogously, one can show by solving the equations
$\sigma^T \rho = 1$ and $\rho \sigma^ T = 1$ inductively in the ring
of lower triangular matrices that $\sigma^T$ is invertible with inverse $\rho$.
\end{proof}

\subsection{A quantum Galois 
group}
Let now $k$ be a commutative ring and $B$ be a $k$-algebra. Let 
$P = \mathrm{End} _{k}(B)$ 
be the ring of $k$-linear maps 
$B \rightarrow B$.
A matrix
$ \sigma \in \Mn( \mathrm{End} _{k}(B))$ 
can alternatively be viewed as
a map 
$B \rightarrow \Mn(B)$, and we 
can demand this to be a ring morphism:

\begin{definition}\label{defautom}
A \emph{quantum
automorphism} of $B$ over  
$k$ 
is a strongly invertible matrix 
$ \sigma \in \Mn(\mathrm{End}_{k}(B))$
satisfying 
\begin{equation}\label{ringmor}
	\sigma_{ij}(1) = \delta _{ij},
	\qquad 
	\sigma_{ij} (ab) = 
	\sum_{l=1}^n \sigma_{il}(a) 
	\sigma _{lj}(b) \qquad 
	\forall a,b \in B.
\end{equation}
Given a quantum automorphism, we denote by
$$
	U_{\sigma} \subset \mathrm{End}_{k}(B)
$$ 
the $k$-algebra generated by the entries 
$ \sigma_{d,ij}$ of the 
$ \sigma _d \in
\Mn(\mathrm{End}_k(B))$,  
$ \sigma _0=\sigma $,
$ \sigma _{d+1} = \hat \sigma_d$. 
\end{definition}  
For $n=1$, this just means that 
$ \sigma $ is a $k$-algebra automorphism of $B$, and that $U_{\sigma}$ is the group
algebra of the subgroup of
the Galois group  
$ \mathrm{Gal}(B/\kk)$ of $B$
over $\kk$ 
that is generated by 
$ \sigma $. In this sense,
the set of quantum automorphisms is a
generalisation of
$\mathrm{Gal}(B/\kk)$. 
 
If $B$ is noncommutative, 
$ \sigma ^{-1}$ is 
in general not a ring
morphism, but note that the set of 
quantum automorphisms is closed 
under $ \sigma \mapsto \hat \sigma $:

\begin{proposition}
If $ \sigma $ is a quantum
automorphism, then so is 
$\hat\sigma$. 
\end{proposition}
\begin{proof}
The key point is to show that 
$\hat \sigma $ is multiplicative. 
To see this, first apply 
$\sigma_{pi}^{-1} $ 
to (\ref{ringmor}) and
sum over $i$. This yields
$$
	\delta _{pj} ab = 	
	\sum_{i,l=1}^n \sigma^{-1} 
	_{pi} (\sigma _{il}(a) \sigma
	_{lj}(b)).
$$
Inserting 
into this equation 
$a= \sigma ^{-1}_{qr}(c),
b=  \sigma^{-1}_{jq}(d)$  
for elements $c,d \in B$
and summing over $j$
and $q$ yields the claim:
\begin{align*}	
	\sum_{q=1}^n 
	\hat\sigma_{rq} (c)  
	\hat \sigma_{qp} (d) 
&= 	
	\sum_{q=1}^n 
	\sigma^{-1}_{qr} (c)  
	\sigma^{-1}_{pq} (d)  = 	
	\sum_{q,j=1}^n 
	\delta_{pj} 
	\sigma^{-1}_{qr} (c)  
	\sigma^{-1}_{jq} (d)  \\
&= 	
	\sum_{i,j,l,q=1}^n \sigma^{-1} 
	_{pi} (\sigma _{il}(
	\sigma^{-1} _{qr} (c)) 
	\sigma_{lj}( 
	\sigma^{-1} _{jq}(d))) \\
&= 	
	\sum_{i,l,q=1}^n \sigma^{-1} 
	_{pi} (\sigma _{il}(
	\sigma^{-1} _{qr} (c)) 
	\delta _{lq}d ) =
	\sum_{i,l=1}^n \sigma^{-1} 
	_{pi} (\sigma _{il}(
	\sigma^{-1} _{lr} (c))d )
 \\
&= 	
	\sum_{i=1}^n \sigma^{-1} 
	_{pi} (\delta _{ir}
	cd)
=	\sigma^{-1} 
	_{pr} (
	cd)
= \hat\sigma_{rp} (cd).\qedhere
\end{align*} 
\end{proof}

However, even when $B$ is
commutative, the (matrix) product 
of quantum
automorphisms is in general not a quantum
automorphism, so quantum automorphisms
do not form groups. As we will
explain in
Section~\ref{interpretation},
they instead generate quantum groups
(Hopf algebras).

\subsection{A quantum subgroup}
Even for basic examples of 
ring extensions, the set of all quantum
automorphisms is usually wild. There are 
various subsets that one can focus on, 
and we will in particular be interested in
the following two attributes:

\begin{definition}
A quantum automorphism $ \sigma $ 
is 
\begin{enumerate}
\item \emph{upper triangular} if 
$
	\sigma _{ij} = 0
$
for $i>j$, and 
\item \emph{locally finite} if 
for all $a \in B$, the set
$
	\{\oX(a) \mid \oX \in U_{\sigma}\} \subseteq B
$ 
is contained in a finitely generated 
$k$-module.
\end{enumerate}
\end{definition}

Both conditions will be
motivated and explained further in
Section~\ref{interpretation}. For now, we
only point out that the upper
triangularity makes it particularly easy
to find such
quantum automorphisms: 

\begin{corollary}\label{upper}
A $k$-algebra morphism
$ \sigma \colon B
\rightarrow \Mn(B)$
with $ \sigma_{ij}=0$ for $i>j$
is a quantum automorphism if and
only if its diagonal entries $\sigma_{ii}$
are invertible for $i = 1, \cdots,n$. In
this case, $U_{\sigma}$ is generated as a $k$-algebra by
the $ \sigma _{ij}$ together with the
$ \sigma _{ii}^{-1}$.
\end{corollary}

\begin{proof}
This follows immediately from Proposition ~\ref{StrInv2}.
\end{proof}

A typical situation in which local
finiteness holds is the following:
 
\begin{proposition}\label{locallyfinite}
Suppose $B$ is a $k$-algebra, and $\{F_d\}_{d \in
\mathbb{Z}}$ is an exhaustive 
filtration
of the algebra $B$ by finitely
generated $k$-modules, $F_d$, 
and that 
$\sigma \in \mathrm{End}_k(B)$ is a quantum automorphism 
with   
$
	\sigma_{ij}(F_d) \subseteq F_d.
$
Then $\sigma$ is locally finite.
\end{proposition}
\begin{proof}
By assumption, all $F_d$ are invariant
under the action of $U_\sigma$, and 
if $b \in B$ is any element, then there
exists $d \in \mathbb{Z} $ with 
$b \in F_d$, and hence $U_\sigma b
\subseteq F_d$. 
\end{proof}

\subsection{Upper triangular quantum automorphisms 
of $
\rational$}\label{rationalfunctionsec}
We will be interested in quantum
automorphisms of coordinate rings of
singular plane curves whose field of
fractional functions is the field 
$k(t)$, and we first classify the upper
triangular quantum automorphisms of the
latter.

\begin{proposition}\label{characterisequk}
For any field $\kk$, the assignment 
$
	\sigma \mapsto 
	\sigma (t) \in \Mn(
	\rational)
$
defines a bijection between 
upper triangular quantum automorphisms 
of $\rational$ over $\kk$  
and upper triangular 
matrices in $\Mn( \rational)$ 
whose diagonal entries are of the 
form 
$$
	\sigma(t) _{ii} = 
	\frac{\alpha _i
	t + \beta _i}{\gamma _i
	t + \delta _i}
$$
for some $\alpha _i,\beta_i,
\gamma _i,\delta _i \in 
\kk$, $ \alpha _i \delta _i - 
\beta _i \gamma _i \neq 0$.
\end{proposition}
\begin{proof}
For any $k$-algebra
$M$, 
$ \sigma \mapsto T:=\sigma (t)$ defines 
a bijection between the set of 
$k$-algebra morphisms 
$ \sigma \colon \polynomial \rightarrow 
M$ and $M$.  
Such an algebra morphism extends 
in at most one way to an algebra morphism 
$ \sigma \colon \rational \rightarrow M$ given
by $ \frac{p}{q} \mapsto p(T) q(T)^{-1}$,  
and it does extend if and only if for any 
$ q \in \polynomial \setminus \!\{0\}$ the element 
$ q(T) \in M$ is invertible in $M$. 

Specialising these general considerations
to the case $M=\Mn(\rational)$, we have
furthermore by elementary linear algebra
over fields:
\begin{enumerate}
\item $p(T)$ is upper triangular for all 
$p \in \polynomial$ if and only if $T$ is so.
\item $q(T)$ is invertible if and only if 
$\det (q(T)) \neq 0$.
\item If $q(T)$ is invertible and upper
triangular, so is $q(T)^{-1}$.
\item If $T$ is upper triangular, 
then $\det (q(T)) = q(T_{11}) \cdots 
q(T_{nn})$. 
\end{enumerate}

We conclude that 
$\kk$-algebra morphisms 
$ \rational \rightarrow 
\Mn(\rational)$ 
that are upper triangular 
correspond bijectively 
to upper triangular 
matrices $T \in \Mn(\rational)$ with 
$ q(T_{ii}) \neq 0$ for all
$i=1,\ldots,n$ and all $q \in \polynomial 
\setminus\!\{0\}$. 

Corollary~\ref{upper} shows 
that such an algebra morphism is a quantum
automorphism if and only if its diagonal
entries $ \sigma _{ii}$ are in the Galois
group of $\rational$ over $\kk$, which is
well-known to be the group of Möbius
transformations \cite[Theorem 7.5.7]{DaCox}.
So if we are given 
an upper triangular matrix $T$ that
defines an upper triangular 
quantum automorphism of $\rational$, 
the $T_{ii}$
are necessarily of the form as stated.  
Conversely, if all $T_{ii}$ are 
of this form, 
then we also have 
$q(T_{ii}) \neq 0$ for all 
$q \in \polynomial \setminus \! \{0\}$,
as the inverse of the unique
$\kk$-algebra automorphism $\sigma _{ii} 
\colon \rational \rightarrow \rational$
that maps $t$ to $T_{ii}$ maps 
$q(T_{ii})$ to $q$; thus $T$ defines an algebra
morphism $\rational \rightarrow \Mn(\rational)$ which
by Corollary~\ref{upper} is a quantum
automorphism. 
\end{proof}

\subsection{Restriction to $\cusp$}
Assume $k$ is a field,  
$ \sigma $ is 
a quantum automorphism 
of $\rational$, and abbreviate as above $T:= \sigma
(t) \in \Mn(\rational)$. Then $ \sigma $
restricts to an intermediate ring 
$\kk \subset B \subset \rational$ 
if and only if for all 
$b = \frac{p}{q} \in B$,
we have
\begin{equation}\label{restrict}
	\sigma (b) = b(T)=p(T) 
	q(T)^{-1} \in
\Mn(B) \subseteq \Mn(\rational).
\end{equation}

The main example we are
interested in is  
$$
	B = \cusp = 
	\mathrm{span}_\kk \{ 
	t^i \mid 
	i \neq 1\} \subset 
	\polynomial,
$$
It is evidently 
sufficient to
test (\ref{restrict}) only for
a set of generators of the
$k$-algebra $B$, so in this
case for 
$b=t^2$ and $b=t^3$. In other words, we
have:
\begin{corollary}
A quantum automorphism 
$ \sigma $ of $\rational$ over $k$ restricts to 
$\cusp$ if and only if 
$ T^2,T^3 \in \Mn(\cusp)$, where 
$ T = \sigma (t)$. 
\end{corollary}

When classifying upper
triangular quantum
automorphisms of $\cusp$
that extend to $\rational$, it
is sufficient to consider
matrices whose entries are
Laurent polynomials:

\begin{proposition}
If an upper triangular matrix 
$T \in \Mn(\rational)$ satisfies 
$T^2,T^3 \in
\Mn(\cusp)$,
then $ T \in \Mn(\laurent)$ whose entries
contain no terms of degree less than 
$-3n+4$.
\end{proposition}
\begin{proof}
We prove  
$T_{ij} \in \laurent$ 
by induction on
$j-i$. 

$j-i=0$: We have shown that the diagonal entries
are Möbius
transformations, i.e.~are of the form 
$T_{ii} = \frac{\alpha _it+ \beta _i}{\gamma _i t +
\delta _i}$. The condition
$T^2 \in \Mn(\cusp)$ and
$\alpha_i \delta_i - \beta_i \gamma_i \neq 0$ force
 $ \beta _i = \gamma _i = 0$, 
thus without loss of generality we have $T_{ii} =  \alpha _i t$.

$j-i = n$: 
Assume that $T_{i\; i+r}$ is for all 
$r < n$ a Laurent
polynomial and contains no term of degree
less than $-3r+1$.

By assumption, the elements
\begin{align*} 
	(T^2)_{i \; n+i} &= 
	(\alpha_i + \alpha_{n+i})t T_{i \; n+i} 
	+ \sum_{r = 1}^{n-1} T_{i \;i+r}T_{i+r
	\; n+i} \\
	(T^3)_{i \; n+i} &= 
	(\alpha_{i}^2 + \alpha_i \alpha_{n+i} + 
	\alpha_{n+i}^2) t^2 T_{i\; n+i} + \\
&\qquad + \sum_{r=1}^{n-1} 
	T_{i \; i+r}(T^2)_{i+r \; n+i} + 
	\alpha_i t T_{i \; i+r}T_{i+r\; n+i} 
\end{align*}
must be in $\cusp$. As 
either $ \alpha _i + \alpha _{n+i} $ or 
$ \alpha _i^2 + \alpha _i \alpha _{n+i} + 
\alpha _{n+i}^2 $ are non-zero, it follows
from the induction hypothesis 
that $T_{i \; n+i}$ is
a Laurent polynomial which contains no
term of degree less than $-3n+1$.
\end{proof}

\begin{remark}
If $k \subsetneq B
\subset \rational$ is any intermediate ring, then
$B$ is a subring of a field, hence an
integral domain, and its fraction field
embeds naturally into $\rational$; so by
Lüroth's theorem, the fraction field is isomorphic
to $\rational$. If $B$ is the coordinate
ring of an algebraic set $V$, this means
that $V$ is an irreducible curve 
which is birationally equivalent
to the affine line. In particular, when
$k$ is algebraically closed, then 
this is the case if and only if $B$ is
finitely generated as a $k$-algebra (by 
Hilbert's Nullstellensatz). 

For $B=\cusp$, the curve $V$ is the cusp 
$$
	V=\{(\alpha , \beta ) \in k^2 \mid 
	\alpha^3 = \beta ^2 \} = 
	\{ (\lambda ^2, \lambda ^3)
\mid \lambda \in k \} \subset k^2,
$$
so in geometric terms, we are
talking about
quantum automorphisms of the cusp that
extend from its coordinate ring to its
field of rational functions. 
\end{remark}

\subsection{Classification 
for $n=2,3$}\label{classificationsection}

Recall that for $l \in \mathbb{N} $
and $ \beta \in \kk$, one defines the
quantum numbers
$$
	\ls l \rs_\beta := 1 + \beta +
	\cdots + \beta^{l-1} =
\frac{1-
\beta^l}{1-\beta},
$$
where the last equality of course only
applies when $ \beta \neq 1$. 

\begin{lemma}\label{n2}
If $z = \sum_{i \in \mathbb{Z}} z_i t^i$,
then the matrix
$$
	T = 
\left(\begin{array}{cc}
	\alpha t & z \\
	0 & \alpha \beta t
\end{array}\right)
$$
corresponds to a quantum
automorphism of $\cusp$ if and only if 
\begin{enumerate}
\item 
$ \ls 2 \rs_ \beta = 0
	\Leftrightarrow \beta =
-1$ and 
$z_{-1}=z_{-3}=z_{-4}=\ldots=0$, or
\item 
$ \ls 3 \rs_ \beta = 0
	\Leftrightarrow \beta =
e^{\pm 2 \pi i/3}$ and 
$z_0=z_{-2}=z_{-3}=\ldots=0$ or
\item $ \ls 2 \rs_\beta, 
\ls 3 \rs_\beta \neq 0$ and 
$z_0=z_{-1}=z_{-2}=\ldots=0$.
\end{enumerate}
\end{lemma}
\begin{proof}
We have
$$
	T^2 = \alpha  
\left(\begin{array}{cc}
	\alpha t^2 & \ls 2 \rs_\beta t z \\
	0 & \alpha \beta^2 t^2
\end{array}\right),\quad
	T^3 = \alpha ^2
\left(\begin{array}{cc}
	\alpha t^3 & \ls  3 \rs_\beta t^2 z \\
	0 & \alpha \beta^3 t^3
\end{array}\right)
$$
so $T$ corresponds to a quantum
automorphism if and only if 
$$
	\ls 2 \rs_\beta tz,\;
	\ls 3 \rs_\beta t^2z \in \cusp
$$
which leads to the conditions as stated. 
\end{proof}

\begin{lemma}\label{n3}
If $z = \sum_{i\in \mathbb{Z} } z_i t^i, \;
y = \sum_{j \in \mathbb{Z}} y_j t^j$, and 
$x = \sum_{l\in \mathbb{Z}} x_l t^l$, 
then the matrix
$$
	T = 
\left(\begin{array}{ccc}
	\alpha t & x & z \\
	0 & \alpha \beta t & y\\
	0 & 0 & \alpha \beta \gamma t
\end{array}\right)
$$
corresponds to a quantum
automorphism of $\cusp$ if and only if 
\begin{enumerate}
\item
\begin{enumerate}
\item
$\beta = -1$ and
$x_{-1}=x_{-3}=x_{-4}=\ldots=0$,
or
\item
$\beta = e^{\pm 2 \pi i/3}$
and
$x_0=x_{-2}=x_{-3}=\ldots=0$,
or
\item 
$ \beta \neq -1,e^{\pm 2
\pi i/3}$ and
$x_0=x_{-1}=x_{-2}=\ldots=0$
\end{enumerate}
\item 
and 
\begin{enumerate}
\item 
$ \gamma =
-1$ and
$y_{-1}=y_{-3}=y_{-4}=\ldots=0$,
or
\item 
$\gamma =e^{\pm 2 \pi i/3}$
and
$y_0=y_{-2}=y_{-3}=\ldots=0$,
or
\item
$ \gamma \neq -1,
e^{\pm 2 \pi i /3} $ and
$y_0=y_{-1}=y_{-2}=\ldots=0$,
\end{enumerate}
\item and 
\begin{enumerate} 
\item
$ \beta \neq - 1- \gamma ^{-1}$
and there are $a,b \in \cusp$
with
$$
	z = (1+\beta + \beta \gamma ) 
	t^{-1} a - t^{-2} b ,\quad
	xy = \alpha (-\ls 3 \rs_{\beta
\gamma } a + \ls 2 \rs_{\beta \gamma
} t^{-1}b),
$$
or
\item
$ \beta = -1-\gamma ^{-1}$
and 
$$	
	c := xy - \alpha \gamma t z
\in \cusp,\ls 3 \rs_ \gamma t c \in \cusp. 
$$
\end{enumerate}
\end{enumerate}
\end{lemma}
\begin{proof}
Applying Lemma~\ref{n2} to
the two matrices 
obtained by deleting from $T$ the first
respectively third
row and column leads to the
conditions (1) and (2). 
Condition (3) 
arises from considering the $(1,3)$-entry
of 
$T^2,T^3$:
$$
	\left(\begin{array}{cc}
	\alpha \ls 2 \rs_{\beta \gamma} &
	1 \\
	\alpha^2 \ls 3 \rs_{\beta \gamma}t   
	&
	\alpha(1+\beta \ls 2 \rs _\gamma)t 
\end{array}\right)
	\left(\begin{array}{c}
	tz \\
	xy 
\end{array}\right) =: 
	\left(\begin{array}{cc}
	a \\
	b
	\end{array}\right)
 	\in \cusp^2
$$
The determinant of the
coefficient matrix is 
$$
	\alpha^2 \beta t ( 
	1 + \gamma + 
	\beta \gamma ).
$$
So we can 
invert the matrix over 
$\kk[t,t^{-1}]$ if
$$
	\beta \neq -1 - \gamma^{-1}.
$$
In this regular case, we obtain
$$
	\frac{1}
	{\alpha^2 \beta ( 
	1 + \gamma + 
	\beta \gamma )}
	 \left(\begin{array}{cc}
	(1+\beta \ls 2 \rs _\gamma)
	t^{-1}  
&
	-t^{-2} \\
	-\alpha \ls 3 \rs_{\beta \gamma}
	&
	\alpha \ls 2 \rs_{\beta \gamma}
	t^{-1} 
\end{array}\right)
	\left(\begin{array}{c}
a \\
b
\end{array}\right) = 
	\left(\begin{array}{c}
	z \\
	xy 
\end{array}\right)
$$
By rescaling $a,b$ this yields
elements $a,b \in \cusp$ with 
$$
	z = (1+\beta + \beta \gamma ) 
	t^{-1} a - t^{-2} b 
$$
and 
$$
	xy = \alpha (- \ls 3 \rs_{\beta
\gamma } a + \ls 2 \rs_{\beta \gamma
} t^{-1}b).
$$
In the singular case
$ 1 + \gamma + \beta \gamma =
0$, the
equation reduces to 
$$	
c: = xy - \alpha \gamma t z
\in \cusp,\quad
	\ls 3 \rs_ \gamma t c \in \cusp. 
$$

The converse is verified by
straightforward computation.
\end{proof}

\subsection{An explicit
example}\label{explicit}
From now on, we assume that
$2,3 \in k$ are invertible 
and that $k$ contains a square
root $i$ of $-1$. 
We will study in detail the
following example of a quantum
automorphism of $\rational$: 

$$ 
	\sigma(t) = T = \begin{pmatrix}
		t & t -i & -\frac{1}{3}t^{-1} - \frac{1}{2}t \\
		0 & -t & t+i \\
		0 & 0 & t
		\end{pmatrix}.
$$ 

This does restrict to 
$B=\cusp$; indeed, we have

$$
	\sigma(x) = T^2 = \begin{pmatrix}
				x & 0 & \frac{1}{3}\\
				0 & x & 0 \\
				0 & 0 & x
			\end{pmatrix}, \; 
	\sigma(y) = T^3 = \begin{pmatrix}
				y & y -ix & -
\frac{1}{2}y \\
				0 & -y & y+ix \\
				0 & 0 & y
			\end{pmatrix},
$$ 
where $x:=t^2, \; y:=t^3$.

By definition, the resulting algebra
$U_{\sigma}$ has four generators that act
as follows on the elements $x,y$:
\begin{gather}
	 \oK := \sigma_{22} : x \mapsto x, \; y \mapsto -y,
	\nonumber\\
	\oE := \sigma _{12} : x \mapsto 0, \; y
	\mapsto y - i x, \qquad	
	\oF := \sigma _{23} : x \mapsto 0, \; y
\mapsto y+ ix,
\label{actongenerators}
\\
	\oZ := \sigma _{13} : x \mapsto \frac{1}{3}, \; y \mapsto - \frac{1}{2}y.
\nonumber
\end{gather}
Since $ \sigma $ is upper
triangular, the operator $\oK$ is an 
algebra automorphism, 
$\oK \in
\mathrm{Gal}(\rational/\kk)$, so 
for all $f,g \in \rational$, we have 
$$
	\oK(fg) = \oK(f) \oK(g),
$$
the
operators $\oE,\oF$ are twisted derivations
satisfying
$$
	\oE(fg) = f \oE(g) + \oE(f) \oK(g),
	\quad \oF(fg) = \oK(f)
\oF(g) + \oF(f)g,
$$
and $\oZ$ is a twisted differential
operator of order 2,
$$
	\oZ(fg) = f \oZ(g) + \oE(f) \oF(g) + \oZ(f) g. 
$$  
This and the action
(\ref{actongenerators}) completely
determines $\oK,\oE,\oF,\oZ$
as $k$-linear maps: 
\begin{lemma}\label{kefzexp}
For all $n \in \mathbb{N} $, we have 
\begin{equation}\label{actionk}
	\oK(t^n) = (-1)^n t^n,
\end{equation}
$$
		\oE(t^n) = \begin{cases}
		t^n - i t^{n-1} &  n \; \mathrm{is \; odd}, \\
		0 &  n \; \mathrm{is \;
even,}
		\end{cases}\quad
		\oF(t^n) = \begin{cases}
		t^n + i t^{n-1} &  n \; \mathrm{is \; odd}, \\
		0 &  n \; \mathrm{is \;
even,}
		\end{cases}
	$$
	$$
		\oZ(t^n) = \begin{cases}
		\frac{n-3}{6}t^{n-2} - \frac{1}{2} t^n &  n \; \mathrm{is \; odd},\\
		\frac{n}{6}t^{n-2} &  n \; \mathrm{is \; even}.
		\end{cases} 
	$$
\end{lemma}
\begin{proof}
This is verified by a straightforward
computation using induction on $n$. 
\end{proof}

In particular, one observes:

\begin{corollary}
The restriction of $ \sigma $
to $B=\cusp$ is locally finite.
\end{corollary}
\begin{proof}
The algebra $k[t^2,t^3]$ inherits a
grading from $k[t]$, so 
$$
	\deg (x) = 2 ,\quad 
	\deg(y) = 3,
$$
and if we denote by 
\begin{equation}\label{deffd}
	F_d:=\mathrm{span}_k
\{1,t^2,t^3,\ldots,t^d\}
\end{equation}
the resulting filtration of $k[t^2,t^3]$,
then by the above formulas 
all assumptions of
Proposition~\ref{locallyfinite} are met.
\end{proof}

\section{The quantum groups 
$H_ \sigma $ and $ A_\sigma $}\label{interpretation}
This section contains the interpretation
and motivation for the above computations:
we explain how quantum automorphisms as
described above give rise to (co)actions
of Hopf algebras on $B$. The
main goal is to prove
that the explicit example discussed in
Section~\ref{explicit} turns
$\cusp$ into a
quantum homogeneous space. 
From now on, we freely use
standard terminology from the theory of
bialgebras and Hopf algebras, see
e.g.~\cite{radford,montgomery,schmudgen,sweedler}
for the necessary definitions. For
simplicity, ``Hopf algebra''
means for us ``Hopf
algebra with bijective
antipode''.

\subsection{The Hopf algebra $H_\sigma$}\label{hsigma}
Recall that a quantum automorphism $
\sigma \in \Mn( 
\mathrm{End}_k(B))  $ of a
$k$-algebra $B$ gives rise
to a Hopf algebra $H_\sigma $ that acts
inner faithfully on $B$, see
\cite{manin,banicabichon,etingofwalton}.
This is constructed as follows:

\begin{enumerate}
\item Consider the free $k$-algebra $k
\langle s_{d,ij} \rangle $ with
generators $s_{d,ij}$, $i,j=1,\ldots,n$, $
d \in \mathbb{Z} $. This carries a unique
bialgebra structure
whose
coproduct and counit are determined by 
\begin{equation}\label{cold}
	\Delta (s_{d,ij}) =
	\sum_{r=1}^n s_{d,ir} \otimes 
	s_{d,rj},\quad
	\varepsilon (s_{d,ij}) = \delta _{ij}. 
\end{equation}
\item Define an action of this 
free algebra on $B$ in which 
the generators act by the
entries of the quantum
automorphisms $ \sigma _d$:
$$
	\act \colon 
	k \langle s_{d,ij} \rangle 
\otimes B \rightarrow B,\quad
	s_{d,ij} \act a :=
	\sigma _{d,ij}(a).
$$
In view of (\ref{cold}), this 
turns $B$ into 
a $k \langle s_{d,ij} \rangle
$-module algebra, that is,
for any $X \in k \langle
s_{d,ij} \rangle $ and $a,b
\in B$, we have 
$$
	X \act (ab) = (X_{(1)} \act
a	)(X_{(2)} \act b).
$$
\item If 
$I \subset k \langle s_{d,ij} \rangle $
is the ideal generated by all elements of
the form 
$$
	\sum_{r=1}^n s_{d,ir} s_{d+1,jr} -
\delta _{ij},\quad 
	\sum_{r=1}^n s_{d+1,ri} s_{d,rj} -
	\delta _{ij}
$$
for some $d,i,j$, then  
$k \langle s_{d,ij} \rangle / I $ becomes
a Hopf algebra with (invertible) antipode
induced by 
$$
	S(s_{d,ij}) := s_{d+1,ji},
$$
and the action $\act$ of 
$ k \langle s_{d,ij} \rangle $ on $B$
descends by
construction to this quotient.
That is, if by abuse of
notation we also denote by 
$ \sigma $ the action
viewed as a morphism 
$$
	\sigma \colon 
	k \langle s_{d,ij} \rangle
\rightarrow 
	\mathrm{End}_k(B),\quad
	X \mapsto X \act -, 
$$
then $I \subseteq
\mathrm{ker}\, \sigma $ so
that $ \sigma $ descends to an algebra
morphism 
$
	k \langle s_{d,ij}  \rangle / I
\rightarrow 
	\mathrm{End}_k(B)
$
that we still denote by 
$ \sigma $.
\item The Hopf 
algebra $ k \langle s_{d,ij} 
\rangle /I$ depends only
on the size $n$ of the matrix 
$ \sigma \in \Mn(
\mathrm{End}_k(B)) $ and not
on the choice of 
$ \sigma $ or $B$.
However, we finally define:
\end{enumerate}
\begin{definition}
We denote by $ H_\sigma $ the
Hopf image  of 
$$ 
	\sigma \colon
	k \langle s_{d,ij}  \rangle / I
\rightarrow 
	\mathrm{End}_k(B).
$$ 
\end{definition} 

That is, $H_\sigma $ is the
universal quotient Hopf algebra 
that acts
on $B$. In other words, $H_\sigma$ is the
quotient of $ k \langle s_{d,ij} \rangle
/I$ by the sum $J_\sigma $ of all Hopf ideals
contained in $ \mathrm{ker}\, \sigma$, see
\cite[Theorem 2.1]{banicabichon} for
further information.

\begin{remark}
More abstractly, 
an algebra morphism 
$ \sigma \colon B \rightarrow \Mn(B)$ is
the same as a measuring of $B$ by 
the coalgebra $C:=\Mn(k)^\star$, the dual of the 
algebra $\Mn(k)$. Steps (1)-(3) construct
the free Hopf algebra with invertible
antipode on this coalgebra, see 
e.g.~\cite{alexandru} and the
references therein. 
The assumption of strong invertibility
implies that the measuring of $B$ by $C$
extends to the free Hopf algebra, and
taking the Hopf image yields the universal
quotient Hopf algebra that has $B$ as a
module algebra. 
\end{remark}

In the upper triangular case,
the results obtained in the
previous section show that 
$H_\sigma$ is a finitely
generated pointed Hopf
algebra:

\begin{proposition}\label{pointedprop}
If $ \sigma $ is upper
triangular, then we have:
\begin{enumerate}
\item $H_\sigma $ is generated
as an algebra by the classes
$[s_{0,ij}]$ with $i \le j$,
together with
$[s_{1,ii}]=[s_{0,ii}]^{-1}$. 
\item $H_\sigma$ is pointed.
\end{enumerate}
\end{proposition}
\begin{proof}
The first claim is shown in the
same way as
Proposition~\ref{StrInv2}. 
The second claim uses a
standard argument: define a
Hopf algebra 
filtration $\{C_f\}$ of $H_\sigma$ by
assigning to $[s_{d,ij}]$ the
filtration degree 
$j-i$, 
$$
	C_f = \mathrm{span}_k 
	\{[s_{d_1,i_1j_1}] \cdots 
	[s_{d_l,i_lj_l}] \mid 
	\sum_q j_q-i_q \le f\}.
$$
This is an algebra filtration
by definition and a coalgebra
filtration as $[s_{d,ij}] =0$ 
if $i>j$. As the $[s_{d,ij}]$
generate $H_\sigma$ as an
algebra, it is exhaustive.  
If $S \subseteq
H_\sigma$ is a simple
coalgebra, then $\dim_k 
S < \infty$, so there exists 
a minimal 
$f \ge 0$ with $S \subseteq 
C_f$, $S \not\subseteq
C_{f-1}$, and if $f>0$, it is
immediately verified that 
$S \cap C_{f-1}$ is a proper
non-zero subcoalgebra of 
$S$, contradicting the fact
that $S$ is simple. 
Finally, if $S \subseteq C_0$
then $S$ is spanned by
group-likes and the span of
any group-like is a
subcoalgebra. So as $S$ is
simple, it is one-dimensional. 
\end{proof}

\subsection{Application to the cusp}\label{applic}
For
the quantum automorphism 
described in
Section~\ref{explicit}, we
abbreviate 
$$
	K:=[s_{0,22}],\quad 
	E:=[s_{0,12}],\quad
	F:=[s_{0,23}],\quad
	Z:=[s_{0,13}] \in H_
\sigma.
$$ 
By
Proposition~\ref{pointedprop}, 
$H_ \sigma $ is generated as
an algebra by these
elements whose images 
in $U_ \sigma $ are
the operators $\oK,\oE,\oF,\oZ$
from Lemma~\ref{kefzexp}. 
Furthermore, the fact that  
$[s_{0,ij}]=0$ for $i>j$
implies that $K$ is
group-like, that is, its
coproduct is given by 
$$	
	\Delta (K) = 
	K \otimes K,
$$
that $E,F$ are $(1,K)$-
respectively $(K,1)$-twisted
primitive,  
$$	
	\Delta (E) = 
	1 \otimes E + 
	E \otimes K,\quad
	\Delta (F) = 
	K \otimes F +
	F \otimes 1,
$$
and $Z$ is of degree 2 with
respect to the coradical
filtration of $H_\sigma$, 
$$
	\Delta (Z) = 1 \otimes Z +
	E \otimes F 
	+Z \otimes 1.
$$

We will now obtain a presentation
of $H_\sigma$ as an algebra.

First, we observe that the
definition of $H_\sigma $
implies:

\begin{lemma}\label{howtofindrels}
We have 
$K^2=1$ and $F=-KE$.
\end{lemma}
\begin{proof}
$K^2-1$ and $KE+F$
are in the kernel of
the representation 
$ \sigma \colon H_\sigma \rightarrow 
U_\sigma$. The coproducts of
these elements are 
$$
	\Delta (K^2-1) = 
	K^2 \otimes K^2 - 1 \otimes 1 = 
	K^2 \otimes (K^2-1) + (K^2-1) \otimes
1,
$$
$$
	\Delta (KE+F) =
	K \otimes (KE+F) +  (KE+F)
\otimes 1,
$$ 
It follows that each one
generates a Hopf ideal
in 
$H_\sigma$ which is in the kernel of 
$ \sigma $, so by definition of $ H_
\sigma $, these elements vanish.   
\end{proof}

Thus $H_\sigma$ is generated
as an algebra by $K,E,Z$. 

Second, we decompose $E$ and
$Z$ into eigenvectors of the
map given by
conjugation by $K$; that is,
we define 
$$
	E_{\pm} := 
	\frac{1}{2}(E \pm
KEK),\qquad
	Z_{\pm} := 
	\frac{1}{2}(Z \pm KZK).
$$
By the definition of these
elements, they 
(anti)commute with $K$:
\begin{lemma}
We have $KE_\pm = \pm E_\pm K$ and 
$KZ_\pm = \pm Z_\pm K$. 
\end{lemma}

The coproduct of these
elements and their action on
$B$ is given by
\begin{gather*}
	\Delta(E_{\pm}) = 
	1 \otimes E_\pm + E_\pm
\otimes K,\\
	\Delta(Z_+) = 
	1 \otimes Z_+ 
	- E_+ \otimes E_+K	
    - E_- \otimes E_-K 
	+ Z_+ \otimes 1,\\
	\Delta(Z_-) = 
	1 \otimes Z_- 
	- E_+ \otimes E_-K	
	- E_- \otimes E_+K 
	+ Z_- \otimes 1,
\end{gather*}
\begin{gather*}
	\sigma (E_+) (t^n) = 
	\begin{cases}
		t^n 
	&  n \; \mathrm{is \; odd}, \\
		0 &  n \; \mathrm{is \;
even,}
		\end{cases}
	\quad
	\sigma (E_-) (t^n) 
	= \begin{cases}
	- i t^{n-1} 
	&  n \; \mathrm{is \; odd}, \\
		0 &  n \; \mathrm{is \;
even,}
		\end{cases}
\\
\sigma (Z_+)(t^n) =\begin{cases}
	\frac{n-3}{6} t^{n-2} - \frac{1}{2}t^n 
	&  n \; \mathrm{is \; odd}, \\
		\frac{n}{6}t^{n-1} &  n \; \mathrm{is \;
even,} 
	\end{cases}
\quad
\sigma (Z_-) (t^n) = \begin{cases}
	0
	&  n \; \mathrm{is \; odd}, \\
		0 &  n \; \mathrm{is \;
even,}
\end{cases}
\end{gather*}

From this, we obtain in a
similar manner as in
Lemma~\ref{howtofindrels}:
\begin{lemma}
We have 
$Z_-=-E_+ E_-$, 
$E_-^2 = 0$, 
$E_+ =- \frac{1}{2}(K-1)$. 
\end{lemma}
\begin{proof}
It follows from the above and
the relation $K^2=1$ that 
the elements $Z_-+E_+E_-$ and 
$E_-^2$ are primitive while
$E_++\frac{1}{2}(K-1)$ is
$(1,K)$-twisted primitive,
and they are straightforwardly
verified to be in $
\mathrm{ker}\, \sigma $, hence
as in
Lemma~\ref{howtofindrels}  
it follows that they vanish in
$H_ \sigma $. 
\end{proof}

So $H_\sigma $ is generated as
an algebra by $K, E_-$ and
$Z_+$. 

Finally, we abbreviate
$$
	Y:= 
	6Z_+ - \frac{3}{2}(K-1), 
	\quad
		D:=iE_-,
	\quad
	 C:= 
	YD - DY. 
$$
Their coproduct is given by
$$
	\Delta(Y) = 
	1 \otimes Y - 
	6D \otimes DK + 
	Y \otimes 1, \quad
	\Delta(C) = 
	1 \otimes C + C \otimes K
$$ 
and they act on $B$ by the operators

\begin{equation}\label{actionyd}
	\oY (t^n) :=\begin{cases}
	(n-3) t^{n-2}, 
		& n \;\mathrm{odd,}\\
	n t^{n-2},  
		& n \;\mathrm{even,}
\end{cases}\qquad
	\oD(t^n):=
\begin{cases}
	t^{n-1} 
	&  n \; \mathrm{is \; odd}, \\
		0 &  n \; \mathrm{is \;
even,}
		\end{cases}
\end{equation}
as well as 
$$
	\oC(t^n) 
:= 
	\begin{cases}
		2 t^{n-3}, 
		& n \;\mathrm{odd,} \\
		0, 
		& n \;\mathrm{even.}
	\end{cases}	
$$

Their commutation relations (as elements in
$H_\sigma$) are as follows:
\begin{lemma}
We have $YK = KY$, $KC = - CK$,
$DC=-CD$,
and 
$$	
	 YC
= CY,\quad C^2=0.
$$ 
\end{lemma}
\begin{proof}
The relations $KY=YK,KC=-CK,DC=-CD$ follow from
the definition of $Y,C,D$ and the
commutation relations already obtained. 
The final 
two relations follow as in
Lemma~\ref{howtofindrels}; 
$YC-CY$ is $(1,K)$-twisted
primitive while $C^2$ is primitive.
\end{proof}

\begin{remark}
Note that we can express
the operators $\oY$ and $\oD$ in terms of $\oK$, 
the differential operator 
$ \frac{d}{dt}$, and the
multiplication operators by 
$t^{-m}$ as 
$$
	\oD = -\frac{1}{2}t^{-1}
(\oK-1),
	\qquad
	\oY=
	t^{-1}\frac{d}{dt}
	+ \frac{3}{2}t^{-2} 
	(\oK-1).
$$ 
The two summands in 
$\oY$ can be
considered separately as
operators on $\rational$ that 
both restrict to $\laurent$, but
only the sum restricts to 
$\cusp$. As operators on 
$\rational$, $\oY_0:=t^{-1}
\frac{d}{dt}$ is a derivation
and $\oY_1:=\frac{3}{2}
t^{-2}(\oK-1)$ is a twisted
derivation,
$$
	\oY_0(fg) =
	f \oY_0(g) + 
	\oY_0(f) g,\quad
	\oY_1(fg) =
	f \oY_1(g) + 
	\oY_1(f) \oK(g),
$$
so it is a rather non-trivial
fact that their sum is a
twisted differential operator
of order 2 on $\cusp$. The
other generator $\oD$ is a twisted
derivation, 
$$
	\oD (fg) = f \oD(g) +
	\oD(f) \oK(g).
$$
\end{remark}

Our aim is to prove that the relations we
have found are complete. In order to do
so, we define the auxiliary Hopf algebra 
\begin{align*}
	\tilde H_\sigma
&:=
	k \langle \tilde K,\tilde D,\tilde Y \rangle / I,\\
	I &:=
	\langle
	\tilde{K^2}-1,
	\tilde K \tilde D + \tilde D \tilde K,\tilde K \tilde Y- \tilde Y \tilde K,
	\tilde {Y^2} \tilde D-2\tilde Y \tilde D \tilde Y+ \tilde D \tilde {Y^2},
	\tilde {D^2}
	\rangle 
\end{align*}
as the algebra generated
by $\tilde K,\tilde D,\tilde Y$ satisfying the relations
established in the lemmata in this
subsection, equipped with the coproduct 
given on generators by the same formulas
as in $H_\sigma$.  
Bergman's diamond lemma \cite{bergman}
immediately yields:
\begin{lemma}\label{basish}
If $\tilde C:=\tilde Y \tilde D-\tilde D\tilde Y$, then the set  
$$
	\{\tilde C^a \tilde D^b \tilde K^c \tilde Y^d \mid  
	a,b,c \in \{0,1\}, d \in 
	\mathbb{N}\} 
$$
is a $k$-vector space basis of $\tilde
H_\sigma$.  
\end{lemma}

We now describe the algebra morphism 
$ \sigma \colon \tilde H_\sigma
\rightarrow U_\sigma $; by showing that
its kernel contains no Hopf ideal, we will
then prove that $\tilde H_\sigma = 
H_\sigma $.

By direct computation, one establishes
that the generators $\oK,\oC,\oD,\oY$ of 
$U_\sigma$ satisfy the following relations
in addition to those satisfied by $\tilde K,\tilde C, \tilde D$
and $\tilde Y$:
\begin{lemma}\label{abovelemma}
We have
$\oC\oD=0$, $\oK\oC=\oC$, $\oK\oD = \oD$.
\end{lemma}

A moment's thought tells that
this a complete presentation
of $U_\sigma$:

\begin{proposition}
The above relations define a
presentation of $U_\sigma$. 
\end{proposition}
\begin{proof}
The claim is that if we define
an abstract algebra 
$k \langle \oK,\oC,\oD,\oY
\rangle /R$, where $R$ is the
ideal generated by the 
relations in
Lemma~\ref{abovelemma} 
together with those that follow from the
ones between $\tilde K,\tilde C,\tilde D,\tilde Y$ in $I$, 
then the resulting
algebra morphism $ k \langle
\oK,\oC,\oD,\oY \rangle /R \rightarrow U_\sigma$
is an isomorphism. 
To do so, observe that using the $k$-vector
space basis of $\tilde H_\sigma$ and the
relations stated in the current
proposition, we obtain a $k$-vector space
basis of
$k \langle \oK, \oC, \oD, \oY \rangle /R$ of the
form 
\begin{equation}\label{basisu}
	\{\oY^a,\oC\oY^b,\oD\oY^c,\oK\oY^d
	\mid
	a,b,c,d \in \mathbb{N} \}.
\end{equation}
It is now straightforward to
show that these operators are
mapped to linearly
independent elements of 
$ \mathrm{End} _k(B)$.
\end{proof}

\begin{remark}
Thus $U_\sigma$ is the Ore
extension of
the 4-dimensional subalgebra
$\mathrm{span}_k\{1,\oC,\oD,\oK\}$
by the derivation 
$\oC \mapsto 0,\oD \mapsto
\oC,\oK \mapsto 0$
(which is given by the commutator with $\oY$). 
\end{remark}
\begin{remark}
Note also that $U_\sigma$
carries a natural grading in
which
$$
	\deg \oK=0,\quad
	\deg \oC=3,\quad
	\deg \oD=1,\quad
	\deg \oY=2,
$$
and that $B$ becomes a graded 
$U_\sigma$-module, 
$(U_\sigma)_i B_j \subseteq
B_{j-i}$.  
\end{remark}

In order to proceed with the
proof that $\tilde H_\sigma 
\cong H_\sigma$, note that by
Lemma~\ref{basish} the subalgebra of
$\tilde H_\sigma $ generated by $\tilde Y$ is 
as an abstract algebra the polynomial
algebra $k[\tilde Y]$ and that $\{\tilde C^a \tilde D^b \tilde K^c \mid 
a,b,c \in \{0,1\}\}$ is a basis of $\tilde
H_\sigma$ as a right $k[\tilde Y]$-module, so as
such, $\tilde H_\sigma $ has rank 8. 
Similarly, $U_\sigma$ becomes a right
$k[\tilde Y]$-module where $\tilde Y$ acts via right
multiplication by $\oY$, and by the
above proposition, $\{1,\oC,\oD,\oK\}$ is
a basis of this $k[\tilde Y]$-module, so this has
rank 4. 
The map $ \sigma \colon \tilde H_\sigma
\rightarrow U_\sigma$ is
right
$k[\tilde Y]$-linear, and we have:

\begin{lemma}\label{basisker}
As a right $k[\tilde Y]$-module, 
$ \mathrm{ker}\, \sigma $ is free with
basis given by 
$$
	\{\tilde C\tilde K+\tilde C,\tilde D\tilde K+\tilde D,
	\tilde C\tilde D,\tilde C\tilde D\tilde K\}.
$$ 
\end{lemma}

As a last ingredient, we
list the one-dimensional
representations of $\tilde
H_\sigma$ and their left and
right
hit actions on $\tilde
H_\sigma$:

\begin{lemma}\label{characters}
For any $ s \in \{-1,1\},\lambda \in k$,
there is an 
algebra morphism
$$
	\chi_{s,\lambda } 
	\colon \tilde H_\sigma 
	\rightarrow k,\quad  
	\tilde K \mapsto s,\quad
	\tilde D \mapsto 0,\quad
	\tilde Y \mapsto \lambda 
$$
and any algebra morphism 
$\tilde H_\sigma \rightarrow
k$ is of this form.
\end{lemma}
\begin{proof}
Immediate. 
\end{proof}

These define
algebra automorphisms 
$\oL_{s,\lambda },
\oR_{s,\lambda } \colon 
\tilde H_\sigma \rightarrow 
\tilde H_\sigma$ given by
$$
	\oR_{s,\lambda } (h) 
	:=
	\chi_{s,\lambda } 
	(h_{(1)}) h_{(2)},
	\quad
 	\oL_{s,\lambda }
	:= 
	h_{(1)} \chi_{s,\lambda } 
	( h_{(2)}). 
$$

On the generators of $\tilde
H_\sigma $, these
automorphisms are given by 
\begin{gather*}
	\oL_{s,\lambda }
	(\tilde K) = 
	\oR_{s,\lambda }
	(\tilde K) = 
	s\tilde K,\\
	\oL_{s,\lambda }
	(\tilde D) = s\tilde D,\quad 
	\oR_{s,\lambda }
	(\tilde D) = 
	\tilde D,\\ 
	\oL_{s,\lambda }
	(\tilde Y) = 
	\oR_{s,\lambda }
	(\tilde Y) = 
	\tilde Y + \lambda.
\end{gather*}
All one-dimensional
representations of $\tilde
H_\sigma$ descend to 
$U_\sigma$:
\begin{lemma}
We have 
$ \mathrm{ker}\, \sigma 
\subseteq \bigcap_{s,\lambda}
\mathrm{ker}\, \chi
_{s,\lambda}$. 
\end{lemma}
\begin{proof}
This follows immediately from 
Lemma~\ref{basisker}. 
\end{proof}
We are now ready to prove:
\begin{proposition}
The quotient $\tilde H_\sigma \rightarrow 
H_\sigma$ is an isomorphism. 
\end{proposition}
\begin{proof}
Assume that $J \subseteq
\mathrm{ker}\, \sigma $ is a
nontrivial Hopf ideal. Then for all 
$h \in J$, we have 
$ \Delta (h) \in J \otimes
\tilde H_\sigma +
\tilde H_\sigma \otimes J$; by the last
lemma, applying 
$ \chi _{s,\lambda }$ to the
left or to the right tensor
component yields an element in
$J$. That is, 
we have 
$$
	\oL_{s,\lambda } 
	(J) = 
	\oR_{s, \lambda } 
	(J) = J.
$$
The maps 
$\oR_{s,\lambda },\oL_{s,\lambda }$
act on the basis elements from 
Lemma~\ref{basish} by 
\begin{align*}
	\oL_{s,\lambda } (\tilde C^a\tilde D^b\tilde K^c\tilde Y^d) 
&=  
	s^{a+b+c} \tilde C^a\tilde D^b\tilde K^c (\tilde Y+\lambda)^d,\\
	\oR_{s,\lambda } (\tilde C^a\tilde D^b\tilde K^c\tilde Y^d) 
&=  
	s^{c} \tilde C^a\tilde D^b\tilde K^c (\tilde Y+\lambda)
^d.
\end{align*}
Thus if 
$$
	X=\sum_{abcd}
	\iota _{abcd}
	\tilde C^a\tilde D^b\tilde K^c\tilde Y^d \in J,
	\quad \iota _{abcd}
	\in k
$$ 
and $d_{\max}(X)$ is the largest 
$d$ such that $ \iota _{abcd} \neq 0$ for
some $a,b,c$, then unless $d_{\max}=0$, 
$$
	X':=X - \oR_{1,1} (X) \in J
$$
is a non-zero element with 
$d_{\max}(X') = d_{\max}(X) -1$. So if $J \neq
0$, it necessarily contains a non-zero element of
the form 
$$
	X=\sum_{abc}
	\iota _{abc}
	\tilde C^a\tilde D^b\tilde K^c.
$$ 
Using now $\oR_{-1,0}$ instead of 
$\oR_{1,1}$, the analogous argument 
shows that $J$ contains a non-zero element of the
form 
$$
	X=\sum_{ab}
	\iota _{ab}
	\tilde C^a\tilde D^b.
$$ 
Considering finally 
$$
	X \pm \oL_{-1,0}(X) 
$$
we find that 
$$
	\iota _{00}
	+ \iota _{11}
	\tilde C \tilde D \in J,\quad
	\iota _{01}
	\tilde C +
	\iota_{10} \tilde D \in J.
$$ 
Since $\oC,\oD \in U_\sigma$
are linearly independent and 
$J \subseteq \mathrm{ker}\,
\sigma$, the second element
vanishes. Using that 
the coproduct of $\tilde C\tilde D$ is 
$$
	\Delta (\tilde C\tilde D) =
	1 \otimes \tilde C\tilde D + 
	\tilde C \otimes \tilde K\tilde D + 
	\tilde D \otimes \tilde C\tilde K + 
	\tilde C\tilde D \otimes 1 
$$
and that $ \tilde C,\tilde D $ are linearly
independent modulo 
$ \mathrm{ker}\, \sigma $, 
one also concludes that the
first element vanishes, a contradiction.
\end{proof}

\begin{remark}
Thus $H_\sigma$ 
does not act faithfully
on $B$ -- for example, 
we now know that 
$CD \in H_\sigma$ is a
non-zero element  while 
$\oC\oD = \sigma (CD) = 0$. 
However, $H_\sigma$ acts by
definition inner faithfully on
$B$, that is, the action 
does descend to algebra, but
not to Hopf
algebra quotients of 
$H_\sigma$. 
\end{remark}

\begin{remark}
Thus $H_\sigma $ is the Ore
extension of the finite-dimensional subalgebra 
generated by 
$K,C,D$ by the derivation given by 
$$
	K \mapsto 0,\quad 
	D \mapsto C,\quad 
	C \mapsto 0. 
$$
In particular, $H_\sigma$ has
Gelfan'd-Kirillov dimension 1, but note
that it is not semiprime (the right ideal
generated by $C$ is a nonzero ideal that
squares to zero) so is not part of the
recent classification of these Hopf
algebras (see
e.g.~\cite{brownzhang,liu} and the
references therein).  
\end{remark}

\begin{remark}
Note that the subalgebra generated by $K$
and $D$ (and similarly the subalgebra
generated by $K$ and $C$) is 
isomorphic as Hopf algebra to Sweedler's
4-dimensional Hopf algebra.  
Recall that 
for any $c \in k$, 
\begin{align*}
	R_c
&:=
	\frac{1}{2} (1 \otimes 1 +
	1 \otimes K + 
	K \otimes 1 - K \otimes K)
\\
&+ \frac{c}{2} (D \otimes D 
	-D \otimes KD + 
	KD \otimes D +
	KD \otimes KD)
\end{align*}
is a universal R-matrix (a
quasitriangular structure) of 
this Hopf subalgebra 
(cf.~\cite[Exercise~12.2.11]{radford}).
This does however not 
define a quasitriangular
structure on $H_\sigma$ as 
$ R_c \Delta (Y) R_c^{-1} \neq 
\Delta^\cop (Y)$, 
but at least for $c=0$,
the corresponding braiding on $B \otimes
B$ is a morphism of $H_
\sigma $-modules. This braiding 
is simply the standard nontrivial
symmetric braiding on the category of
graded vector spaces,
$$
	t^i \otimes t^j \mapsto
	(-1)^{ij} t^j \otimes t^i.
$$
\end{remark}

\subsection{The Hopf algebra $A_ \sigma $}
We now pass to a dual picture: 
assume that $ \sigma $ is a locally finite quantum
automorphism of a $k$-algebra $B$.
Then the $H _\sigma $-action on $B$ arises
by dualisation from an $H_ \sigma^\circ
$-coaction, where 
$H_ \sigma ^\circ $ denotes the Hopf dual
of $H_ \sigma $, that is, the universal Hopf
algebra contained in the vector space dual 
$ H_ \sigma ^\star$. In other words, $B$ is a
right $H_\sigma^\circ $-comodule algebra
with a coaction that we denote by 
$$
	\rho \colon B \rightarrow B \otimes H_ \sigma^\circ
,\quad b \mapsto b_{(0)} \otimes b_{(1)}.
$$ 

\begin{definition}
We denote by $A_\sigma \subseteq
H_\sigma^\circ $ the Hopf subalgebra
generated by the matrix coefficients 
$\{f(b_{(0)}) b_{(1)} \mid 
b \in B, f \in B^\star\}$ of $ \rho $.  
\end{definition}


When $ H_ \sigma $ is
infinite-dimensional, $A_\sigma $ could be a
proper Hopf subalgebra of $ H_ \sigma
^\circ $, but note that it is always
dense:
\begin{proposition}
The restriction of the pairing of $H_
\sigma ^\circ \otimes H_ \sigma \longrightarrow k $ to 
$A_\sigma \otimes H_ \sigma $ is non-degenerate.
\end{proposition}
\begin{proof}
The degeneration space 
$$
	\{ X \in H_ \sigma \mid 
	a(X) =0 \; \forall \; a \in A_\sigma\}
$$
is a Hopf ideal that acts trivially on
$B$, hence vanishes by the definition of
$H _\sigma $. 
\end{proof}

Note that if $M \subseteq B$ is any
finite-dimensional $H_\sigma$-submodule
that generates $B$ as an algebra, then we have 
$$
	B \cong TM / R 
$$
as an $H _ \sigma $-module algebra, where
$TM$ is the tensor algebra of $M$ (over
$k$) and $R$ is the 2-sided ideal of
relations that hold among the elements of
$M$ in the algebra $B$. Then we have:

\begin{lemma}\label{asisgen}
$A_\sigma$ is
generated as a Hopf algebra by the matrix
coefficients of $M$. 
\end{lemma}
\begin{proof}
The matrix coefficients of $M^{\otimes n}$ 
are sums of products of $n$ matrix
coefficients of $M$, and the space of
matrix coefficients of a quotient comodule 
$M^{\otimes n} /(R \cap M^{\otimes n})$ is
a subspace of the space of matrix
coefficients of $M^{\otimes n}$. 
\end{proof}

Alternatively, this can be
formulated in coordinates:
if  
$e_1,\ldots,e_{\dim_k M}$ is a vector
space basis of $M$, then $A_ \sigma $ is
generated as a Hopf algebra by the
functionals
$a_{ij} \in H_\sigma^\circ $,
$i,j=1,\ldots,\dim_k M$, for which 
$$
	X \act e_i=\sigma (X) (e_i) =	
	\sum_{j=1}^{\dim_k M} a_{ji}(X) e_j,\quad
	X \in H_ \sigma,
$$ 
so the coaction $ \rho $ is given by
$$
	M \rightarrow M \otimes A_\sigma,\quad
	e_i \mapsto \sum_{j=1}^{\dim_k M} e_j \otimes
	a_{ji}. 
$$

\begin{remark}\label{pointeda}
In general, the matrix coefficients
of $M$ do not generate $A_\sigma$ as an
algebra -- the subalgebra that they
generate is a subbialgebra of $A_\sigma$
as the span of the matrix coefficients is
a subcoalgebra, but this subbialgebra is
not closed under the antipode in general.
However, if the matrix coefficients can be
chosen to  be upper triangular (that
is, there is a vector space basis of $M$
such that $a_{ij}=0$ for
$i>j$), then the same
arguments that were used in
Proposition~\ref{pointedprop}
show that $A_\sigma $ is generated by
the 
$a_{ij}$ together with the $a_{ii}^{-1}$,
and that 
$A_\sigma $ is a pointed Hopf algebra.
\end{remark}

\subsection{Embedding $B$ into 
$A_ \sigma $}
Now assume that $ \chi \colon B
\rightarrow k$ is an algebra map, that is,
a one-dimensional
representation of $B$.
This induces a map (see
\cite[Section
3]{etingofwalton} for a more
detailed 
discussion of this map) 
$$
	\iota 
	:=(\chi \otimes 
	\mathrm{id} _{A_\sigma})
\circ \rho 
\colon 
	B \rightarrow A_\sigma,\quad 
	b \mapsto 
	\chi (b_{(0)}) b_{(1)}.   
$$
By definition, this is a morphism of
algebras and of right
$A_\sigma$-comodules and
hence it maps $B$ to a right coideal
subalgebra of
$A_\sigma$.

\begin{proposition}
$ \iota $ is injective if and only if for
all $b \in B$, $ b \neq 0$, there exists 
$X \in H _\sigma$ with 
$ \chi (X \act b) \neq 0$.  
\end{proposition}
\begin{proof}
The map $ \iota $ is not
injective if there exists $b
\in B$, $b 
\neq 0$,
with 
$$
	\iota (b) = \chi (b_{(0)}) 
	b_{(1)} = 0.
$$
This is an element in $A_\sigma \subseteq
H_\sigma
^\circ$, so it is zero if and only if it
pairs trivially
with all elements $ X \in H_\sigma$. Thus 
$ \iota $ is not injective if and only if there
exists $b \in B$, $ b \neq 0$, such that 
$$
	X(\iota (b)) = \chi (b_{(0)}) 
	b_{(1)}(X) = \chi ( X \act b)=0
$$
for all $X \in H_\sigma$. 
\end{proof}

If this condition is satisfied, then $
\iota $ embeds $B$ as a right
coideal subalgebra into $A_
\sigma$. 
In particular, when $B=k[V]$ is the
coordinate ring of an algebraic set $V$,
then $ \chi $ corresponds to a point $p
\in V$. The above proposition states
that $B$ can be embedded as a right
coideal subalgebra into 
$A_\sigma$ provided that there exists a
point $p \in V$ such that for any non-zero
regular function $b \colon V \rightarrow
k$ there exists some $X \in H_\sigma $
such that the function $X \act b$ does not
vanish at $p$. 

\subsection{The cusp again}\label{applic2}
To compute a full presentation of
$A_\sigma$ in a given example is tedious,
but relatively
straightforward. Like
elsewhere, we illustrate the
theory with our main example:

\begin{proposition}\label{propas}
If $B=\cusp$ and $M$ is the 
$H_\sigma$-module $F_3$ from 
$(\ref{deffd})$
with basis 
$e_1=1,e_2=t^2,e_3=t^3$,
then we have:
\begin{enumerate}
\item There is a surjective
algebra morphism  
$$
	\pi \colon 
	k \langle \gamma, \varphi ,
	\psi \rangle
	\rightarrow A_\sigma
$$ 
given by $ \pi (\gamma) = 
a_{13}$, $ \pi (\varphi) = 
a_{23}$, $ \pi (\psi) =
a_{33}$ whose kernel 
is the ideal
generated by
\reqnomode
\begin{gather} 
\begin{aligned}\label{relationsas}
	 \psi^2 - 1, \quad 
	\gamma \psi + \psi \gamma
,\quad
	\varphi \psi + \psi \varphi , \\
	27 \gamma^2 - \varphi^6, \quad 
	3(\gamma \varphi + \varphi
\gamma) -
	\varphi^4.
\end{aligned}
\end{gather}
\item In this presentation, 
the coalgebra structure of
$A_\sigma$ is given by
\begin{gather}\label{coalgebraas}
	\Delta(\psi) = \psi \otimes \psi, 
	\quad 
	 \Delta(\varphi) = 
	1 \otimes \varphi + 
	\varphi \otimes
\psi,\nonumber\\ 
\Delta(\gamma) = 
	1 \otimes \gamma + 
	\frac{1}{3}\varphi^2 \otimes \varphi + 
	\gamma \otimes \psi,\\
	\varepsilon(\psi) = 1, \quad 
	\varepsilon(\gamma) = 
	\varepsilon(\varphi) = 0,
\nonumber
\end{gather}
and its antipode is given by 
\reqnomode
\begin{align} \label{antipodeas}
	S(\psi) = \psi,\quad
	 S(\varphi) = -\varphi \psi , \quad 
	S(\gamma) = ( 
	\frac{1}{3}\varphi^3
- \gamma)
	\psi.
\end{align}
\end{enumerate}
\end{proposition}

For the proof that will be
split into several lemmata, we 
introduce a
redundant generator 
$\delta$ and first observe:

\begin{lemma}
Let 
$J \lhd k \langle \gamma ,
	\varphi , \psi , \delta
	\rangle$ be the ideal
generated by the elements
$(\ref{relationsas})$ together  
with $ \delta - \frac{1}{3}
\varphi ^2$. Then we have:
\begin{enumerate}
\item $k \langle \gamma ,
\varphi , \psi , \delta
\rangle $ carries a unique
bialgebra structure such that
\reqnomode
\begin{gather}
	\Delta(\psi) 
	= \psi \otimes \psi, 
	\quad 
	 \Delta(\varphi) = 
	1 \otimes \varphi + 
	\varphi \otimes
	\psi,\nonumber\\ 
	\Delta(\gamma) = 
	1 \otimes \gamma + 
	\delta \otimes \varphi + 
	\gamma \otimes \psi,\quad
	\Delta (\delta ) = 
	1 \otimes \delta + \delta
	\otimes 1,\\
	\varepsilon(\psi) = 1, \quad 
	\varepsilon(\gamma) = 
	\varepsilon(\varphi) = 
	\varepsilon (\delta ) = 0.
\nonumber
\end{gather}
\item The ideal $J$ is a coideal, so this
bialgebra structure descends to 
$$
	\tilde A_\sigma := 
	k \langle \gamma, \varphi ,
	\psi , \delta \rangle /J.
$$
\item The bialgebra $\tilde A_\sigma $ is
a Hopf algebra whose antipode is given on
the generators by $(\ref{antipodeas})$ and 
$S( \delta ) = - \delta $.  
\end{enumerate}
\end{lemma}
\begin{proof}
All claims are verified by straightforward
computations that we leave to the reader.
\end{proof}

Next, we note:
\begin{lemma}
There is a surjective bialgebra morphism 
$$
	\pi \colon k \langle \gamma , \varphi ,
\psi , \delta \rangle \rightarrow A_\sigma
$$	
satisfying 
$ \pi (\gamma ) = a_{13}$, 
$ \pi ( \varphi ) = a_{23}$, 
$ \pi (\psi ) = a_{33}$, and 
$ \pi ( \delta ) = a_{12}$. 
\end{lemma}
\begin{proof}
Recall first that the values of the
functionals $ a_{ji} $ on the generators
$K,D,Y$ of $H_ \sigma $ are by
(\ref{actionyd}) and (\ref{actionk}) given
by the following matrices (for later use,
we also list the values on $C$):
\reqnomode
\begin{gather}\label{matrixpairing}
\begin{aligned}
	a(K) &= 
	\begin{pmatrix}
	1 & 0 & 0 \\
	0 & 1 & 0\\
	0 & 0 & -1
	\end{pmatrix},\quad
	& a(D) &= 
	\begin{pmatrix}
	0 & 0 & 0 \\
	0 & 0 & 1\\
	0 & 0 & 0
	\end{pmatrix},\\
	a(Y) &= 
	\begin{pmatrix}
	0 & 2  & 0 \\
	0 & 0 & 0\\
	0 & 0 & 0
	\end{pmatrix},\quad
	& a(C) &= 
	\begin{pmatrix}
	0 & 0  & 2 \\
	0 & 0 & 0\\
	0 & 0 & 0
	\end{pmatrix}.
\end{aligned}
\end{gather}
As the generators of $H_\sigma$ act by upper
triangular matrices, all elements in $H_
\sigma $ act by upper triangular matrices,
hence as elements of 
$H_\sigma ^\star$, the $a_{ij}$ with $i>j$ 
vanish (this is just a
restatement of the fact that $F_i$ is a
filtration of $B$ by $H_ \sigma
$-submodules). 

Finally, $a _{11} = a _{22} =
\varepsilon_{H_\sigma} = 1_{A_\sigma}$ is
the counit of the coalgebra $H_\sigma$ hence the unit of
the algebra $A_\sigma$. Using this it is
immediately verified that
the algebra morphism $ \pi $ 
defined on the generators $ \gamma ,
\varphi , \psi , \delta $ as in the lemma
is a bialgebra morphism, and 
Lemma~\ref{asisgen} implies that $ \pi $
is surjective. 
\end{proof}

Next, we show that $ \pi $ descends to
a Hopf algebra surjection 
$ \tilde A_\sigma \rightarrow A_\sigma$:

\begin{lemma}
$ J \subseteq \mathrm{ker}\, \pi
$. 
\end{lemma} 
\begin{proof}
The bialgebra map $ \pi $ induces a
pairing of bialgebras 
$$
	\langle -,- \rangle \colon 
	H_\sigma \otimes 
	k 
	\langle \gamma , \varphi ,
\psi , \delta \rangle  \rightarrow k,
\quad
	\langle X,a \rangle :=
	\pi (a) (X).
$$
To prove the lemma, one has to show that
this pairing descends to a
pairing between $H_\sigma $ and  
$\tilde A_\sigma$. 
In order to do so, it is
sufficient to show that for 
each of the six relators $ \xi $ that
generate $J$ and for all $i,j,k \in
\{0,1\}, l \in \mathbb{N} $, we have 
$$
	\langle C^iD^jK^kY^l, \xi \rangle =0.
$$ 
This is verified by  
straightforward computation.  
For example, using that $ \langle -,-
\rangle $ is a pairing of bialgebras, one
obtains 
\begin{align*}
	\langle C^i D^j K^k Y^l , \varphi^2
\rangle 
&=
	\langle (C^iD^jK^kY^l)_{(1)} , 
	\varphi 
\rangle 	
	\langle (C^iD^jK^kY^l)_{(2)} , \varphi 
\rangle 	\\
&=
	\langle C^i_{(1)} D^j_{(1)} K^k_{(1)} Y^l_{(1)} , 
	\varphi 
\rangle 	
	\langle C^i_{(2)} 
	D^j_{(2)} K^k_{(2)} Y^l_{(2)} , \varphi 
\rangle 	\\
&=
	(\langle C^i_{(1)} , 
	\varphi _{(1)}
	 \rangle 
	\langle  D^j_{(1)} , \varphi _{(2)}
	 \rangle 
	\langle K^k_{(1)} , 
	\varphi _{(3)}
	\rangle 
	\langle Y^l_{(1)} ,\varphi _{(4)}
	\rangle )\\
& \quad	(\langle C^i_{(2)} , \varphi _{(1)}
 \rangle 
	\langle  D^j_{(2)} , \varphi _{(2)}
 \rangle 
	\langle K^k_{(2)} , \varphi _{(3)}
 \rangle 
	\langle Y^l_{(2)} , \varphi _{(4)}
	\rangle).
\end{align*}
Inserting the explicit coproducts and at
last the values 
(\ref{matrixpairing}) 
of the pairings of the generators, one 
obtains that the above 
is equal to 
$$
	\langle C^i D^j K^k Y^l , 3 \delta
\rangle =
	3 \delta_{i0}
	\delta _{j0} 
	\delta _{l1}
$$
so that $ \langle - , \delta - \frac{1}{3}\varphi ^2
\rangle $ vanishes as a $k$-linear
functional on $H_\sigma$. The other five
relators are treated in the
same way. 
\end{proof}

The proof that  
$ \pi \colon \tilde A_\sigma \rightarrow
A_ \sigma $ is also injective relies on:

\begin{lemma}
The set 
$$
	\{ \gamma ^a \varphi ^b \psi ^c \mid 
	a,c \in \{0,1\}, b \in \mathbb{N} \}
$$
is a $k$-vector space basis of $\tilde
A_\sigma$. 
\end{lemma}

\begin{proof}
Like Lemma~\ref{basish}, this is a
standard application of Bergman's diamond lemma.
\end{proof}

Thus to prove the injectivity 
of $ \pi $, one has to show that 
the elements $ \pi ( \gamma ^a \varphi ^b
\psi ^c ) \in A_\sigma$ are linearly
independent over $k$. This is maybe shown
most easily by explicitly computing the
values of the functionals:
\begin{lemma}
The dual pairing  
$
	\langle -,- \rangle \colon
	H_\sigma \otimes \tilde A_\sigma 
	\rightarrow k
$
satisfies
$$
	\langle C^i D^j K^k Y^l , 
	\gamma ^a \varphi ^b \psi ^c \rangle =
	\delta _{j+2l,b} 
	\delta _{i a} 
	(-1)^{j (a+c) + ic+
	ab+k (a+b+c)}
	2^a 
	6^l l!
$$
\end{lemma}
\begin{proof}
A direct computation and a nested induction on
$b$ and $l$ shows 
\begin{align*} 
	\Delta ( \gamma ^a \varphi ^b \psi ^c) 
&= 
	(1 \otimes \gamma + 
	\frac{1}{3} \varphi ^2 \otimes \varphi 
	+ \gamma \otimes \psi)^a \\
& \quad \cdot 
	\sum_{l=0}^b 
	(1 - \ls  (b+1)l \rs_{-1})
\begin{pmatrix}
	[b/2] \\
	[l/2]
	\end{pmatrix}
	\varphi ^l \psi ^c
	\otimes 
	\varphi ^{b - l}
	\psi^{\ls bl \rs_{-1}+c},
\end{align*}
where as before 
$\ls n \rs_q = 1 + q + \cdots + q^{n-1}$,
which for $q=-1$ is $0$ if $n$
is even and $1$ if $n$ is odd, and 
$$
	[n/2]:=
	\begin{cases} 
	(n-1)/2 & n \mbox{ odd},\\
	n/2 & n \mbox{ even.}
	\end{cases} 
$$ 

Using this and the formulas for the
coproduct of and the relations between the
generators of $H_\sigma$ respectively 
$ \tilde A_\sigma$ as well as their
pairing (\ref{matrixpairing}), one computes
\begin{align*}
 \quad\langle C^i D^j  Y^l K^k, \gamma ^a
	\varphi ^b \psi ^c \rangle 
&=	\langle C^i D^j Y^l \otimes K^k , 
	\Delta (\gamma^a \varphi ^b \psi ^c)
\rangle \\
&= \langle C^iD^jY^l, 
	\gamma ^a\varphi ^b \psi ^c \rangle 
	\langle K^k ,
	\psi ^a \psi^{\ls b^2 \rs_{-1}+c}
\rangle \\
&=
\langle C^i \otimes D^jY^l, 
	\Delta (\gamma ^a\varphi ^b
\psi ^c) \rangle 
	(-1)^{k (a+c+ \ls b^2 \rs_{-1})}. 
\end{align*}
If $i=0$, the above is equal to
\begin{align*} 	
\ldots &= \langle D^jY^l, 
	\gamma ^a\varphi ^b
\psi ^c \rangle 
	(-1)^{k (a+c+ \ls b^2
\rs_{-1})} \\
	&=\langle \Delta (D^jY^l) , 
	\gamma ^a \otimes \varphi ^b
\psi ^c \rangle 
	(-1)^{k (a+b+c)} \\
&= 
	\delta _{ia} 
	\langle D^jY^l , 
	\varphi ^b \psi ^c \rangle 
	(-1)^{k (a+b+c)}.
\end{align*}
The last equality follows since the pairing of the
first tensor component of 
$ \Delta (D^j Y^l)$ with  
$ \gamma $ vanishes.

If $i=1$, we have instead
\begin{align*}
 \ldots &=	
	\langle C \otimes D^jY^l, 
	\Delta (\gamma ^a\varphi ^b
\psi ^c) \rangle 
	(-1)^{k (a+c+ \ls b^2
\rs_{-1})} \\
&= \langle D^jY^l, 
	\varphi ^b
\psi ^{a+c} \rangle 
	2\delta _{i a} (-1)^{ic}
(-1)^{ac}
	(-1)^{k (a+c+ \ls b^2
\rs_{-1})} \\
 &= \langle D^jY^l, 
	\varphi ^b
\psi ^{a+c} \rangle 
	2\delta _{i a} (-1)^{ic+
	ab+k (a+b+c)}.
\end{align*}
We have deliberately written $i$ instead of $0$
respectively $1$ in the above two cases as
we now can merge them again:
\begingroup
\allowdisplaybreaks
\begin{align*}
	\langle C^i D^j  Y^l K^k, \gamma ^a
\varphi ^b \psi ^c \rangle 
 &= \langle D^jY^l, 
	\varphi ^b
\psi ^{a+c} \rangle 
	2^a\delta _{i a} (-1)^{ic+
	ab+k (a+b+c)}\\
 &= \langle \Delta (D^jY^l), 
	\varphi ^c \otimes 
\psi ^{a+c} \rangle 
	2^a\delta _{i a} (-1)^{ic+
	ab+k (a+b+c)}\\
 &= \langle D^j \otimes Y^l, 
	\Delta (\varphi ^b) \rangle
	2^a  
	\delta _{i a} 
	(-1)^{j (a+c) + ic+
	ab+k (a+b+c)}\\
 &= 6^l l! (1- \ls j+2l+1\rs_{-1})^j
	2^a \delta _{j+2l,b} 
	\delta _{i a} \\
&\quad \cdot  
	(-1)^{j (a+c) + ic+
	ab+k (a+b+c)} \\
& =  \delta _{j+2l,b} 
	\delta _{i a} 
	(-1)^{j (a+c) + ic+
	ab+k (a+b+c)} 6^l l! 
	2^a.\qedhere
\end{align*}
\endgroup
\end{proof}

Therefore, if we define 
$$
	E_{uvw} := \frac{(-1)^{2[v/2](u+v) -v -uw-uv}}{6^{[v/2]} [v/2]! 2^{u+1}}
	C^u D^{v-2[v/2]} Y^{[v/2]} 
	(1 + (-1)^{u+v+w} K), 
$$
then these elements also form a basis of
$H_\sigma$, and we have
$$
	\langle E_{uvw} , \gamma ^a \varphi ^b
\psi ^c \rangle = 
	\delta _{ua} \delta _{vb} \delta _{wc}.
$$
This shows that the elements 
$ \pi ( \gamma ^a
\varphi ^b \psi ^c) \in A_ \sigma $ are
linearly independent which finishes the proof of
Proposition~\ref{propas}. 

\begin{remark}
In our example, $ \sigma $ is a $3 \times
3$-matrix with entries in $
\mathrm{End}_k(B)$ and $(a_{ji})$ is a 
$3 \times 3$-matrix with entries in 
$H_\sigma^\circ $, but be aware it is a pure
coincidence that these sizes match.
\end{remark}

Finally, we prove that the map 
$ \iota \colon B \rightarrow A_\sigma$ is
injective for every point $p$ on the cusp:
the algebra morphisms $ \chi \colon 
B=k[t^2,t^3] \rightarrow k$ are in bijection with 
the points $p=(\lambda ^2,\lambda ^3) $ on
the cusp, $ \lambda \in k$, and
the algebra morphism $ \iota $ is given on
generators by
$$
	t^2 \mapsto 
	\lambda^2 1 + \frac{1}{3} \varphi ^2,
	\quad 
	t^3 \mapsto 
	\gamma +
	\lambda^2 \varphi + 
	\lambda^3 \psi. 
$$

Now we observe:

\begin{lemma}
The map $ \iota \colon B \rightarrow
A_\sigma $ is injective for all algebra
maps $ \chi \colon B \rightarrow k$.
\end{lemma}

\begin{proof}
The element 
$X:=Y+D \in H_\sigma$
acts by 	
$$
		\oX (t^n) := \begin{cases}
		(n-3) t^{n-2} + t^{n-1}, 
	&  n \; \mathrm{is \; odd},\\
		nt^{n-2} & 
	n \; \mathrm{is \;
even}.
		\end{cases}
	$$
In particular, if $b \in F_d
\setminus F_{d-1}$ is a
polynomial of degree $d>0$, then 
applying $\oX^{d/2}$ 
when $d$ is even and
$\oX^{(d+1)/2}$
when $d$ is odd yields a
non-zero scalar, that is, 
a 
constant regular
function on the cusp. 
\end{proof}

\begin{remark}\label{greatnewgenerators}
To clean up the presentation
of $A_ \sigma $, we finally
introduce a new set of
generators:
$$
	\oa :=
	\frac{1}{6}\varphi, \quad 
	\ob := \frac{1}{6^2}\gamma -
\frac{1}{6^3}\varphi^3,\quad
	\oc := \psi.
$$
The new relations between the
new generators are
\begin{gather*}
	\oa \ob + \ob \oa =  
	\oa \oc + \oc \oa = 
	\ob \oc + 
	\oc \ob = 0,\\
	\quad 3 \ob^2 =
	\oa^6, \quad \oc^2=1.
\end{gather*}
Their coproduct, counit, and antipode are
\begin{gather*}
	\Delta(\oa) = 1 \otimes \oa
+ \oa \otimes \oc,\quad
	\Delta (\oc) = \oc \otimes 
	\oc \\
	\Delta(\ob) = 1 \otimes \ob + 
	\oa^2 \otimes \oa - \oa
\otimes \oa^2 \oc + \ob
\otimes \oc \\
	\varepsilon(\oa) = 
	\varepsilon(\ob) = 0, \quad
	\varepsilon (\oc) = 1,\\
	S(\oa) = - \oa \oc, \quad 
	S(\ob) = - \ob \oc,\quad
	S(\oc) = \oc.
\end{gather*}
It is now evident that $A_
\sigma $ is a graded Hopf
algebra with $\oa,\ob,\oc$ of
degree $1,3,0$, respectively.
The monomials
$\oa^l\ob^m\oc^n$, $l \ge 0$, 
$m,n \in \{0,1\}$ are a vector
space basis, so $A_ \sigma $
has GK-dimension 1. One also
observes that $Z(A_\sigma )$
is the polynomial ring 
$\kk[\oa^2]$. 
In terms of the new
generators, the embedding of
the cusp is given by
$$
	t^2 \mapsto 
	\lambda^2 1 + 
	12 \oa^2,
	\quad 
	t^3 \mapsto 
	6\lambda^2 \oa +
	36 \oa^3 +  	
	36 \ob +
	\lambda^3 \oc. 
$$
\end{remark}

\subsection{Faithful flatness}
To complete the proof that
the embedding $ \iota \colon 
B \rightarrow A_ \sigma $
turns $B$ into a quantum
homogeneous space
in the sense of 
\cite{mullerschneider},
we remark that a
result of Masuoka implies that 
$A_\sigma$ is a faithfully flat 
$B$-module. 

As a preliminary result, we
need to compute the coradical
of $A_\sigma$:

\begin{proposition}
The Hopf algebra $A_\sigma$ is
pointed. Furthermore,  
for any point $p$ in the cusp, $A_\sigma$
is faithfully flat over 
$ \iota (B)$.
\end{proposition}
\begin{proof}
We just observed in
Remark~\ref{greatnewgenerators} 
that $A_ \sigma $ is
a graded Hopf algebra with a
vector space basis given by the
monomials $\oa^l\ob^m\oc^n$ 
which are of degree $l+3m$. It
follows that the degree 0
part of $A_ \sigma $ is the group
algebra
$\mathrm{span}_\kk\{1,\oc\}
\cong \kk \mathbb{Z} _2$.  

By
\cite[Theorem (1)]{masuokaCoco}, 
we only need to prove that the
intersection of $ \iota (B)$
with the coradical 
$\mathrm{span}_k\{1, \oc \}$
is invariant under the
antipode $S$ of $A_\sigma$. 
However, the
restriction of the antipode
$S$ to the coradical is the identity map
($ 1 $ and $ \oc $
are their own inverses), so
there is nothing to prove. 
\end{proof}

\subsection{$*$-structures and
involutions}\label{starstructures}
As we have seen in Lemma~\ref{n3},  
even the classification of 
the $3 \times 3$-upper triangular quantum automorphisms
of $B = \cusp$ is rather involved. The
explicit example studied since
Section~\ref{explicit} was obtained by 
demanding in addition that the quantum
automorphism is compatible with a chosen
$*$-structure on $B$. This 
might be of interest in its own
right, for example as it is the starting
point for the transition from
the algebraic theory of Hopf
algebras acting on rings to
the analytic theory of 
locally compact
quantum groups acting on
$C^*$-algebras. 

\begin{definition}
Assume that $k$ is a field with a chosen
field automorphism 
$k \rightarrow k,\lambda \mapsto \bar \lambda
$ that is involutive, i.e.~which satisfies 
$\bar {\bar \lambda } = \lambda $. 
\begin{enumerate}
\item For each $k$-vector space $V$, we denote by 
$\overline{V}$ the \emph{conjugate} vector space
which is the same abelian group but whose
scalar multiplication is twisted by 
$\bar \cdot$ to $ \lambda
\cdot_{\overline{V}} v
:= \bar \lambda v$, $\lambda \in k,v \in
V$.
\item A \emph{$*$-structure} on a
$k$-algebra $B$ is an involutive 
$k$-algebra isomorphism $ * \colon B
\rightarrow \bar B^\op$.
\item An \emph{involution} on a
$k$-algebra $P$ is an involutive
$k$-algebra isomorphism
$\theta \colon P \rightarrow \bar P$. 
\end{enumerate}
\end{definition}
So explicitly, a $*$-structure is a map 
$B \rightarrow B$, 
$b \mapsto b^*$ satisfying 
for all $a,b \in B, \lambda \in k$
$$
	(\lambda a + b)^* = \bar \lambda a^*+ b^*,\quad
	(ab)^* = b^*a^*, \quad 
	a^{**} = a,
$$
while an involution is a map $\theta \colon P
\rightarrow P$ such that 
for all $ \lambda \in k, f,g
\in P$ 
\begin{gather*}
	\theta (\lambda f + g) = 
	\bar \lambda \theta(f) +
	\theta(g),\quad
	\theta(f g) =
	\theta(f)  
	\theta(g),\quad
	\theta(\theta(f)) = f.
\end{gather*}

Typically, these notions are 
considered for $k=
\mathbb{C} $ with involution
given by complex conjugation, 
see
e.g.~\cite[Section~1.2.7]{schmudgen}.

The following is immediate:
\begin{lemma}
A $*$-structure on an algebra $B$ induces
an involution on 
$ P=\mathrm{End}_k(B)$ given by
${\thetab}(f) := * \circ f \circ *$.
\end{lemma}

On Hopf algebras, one demands
the following compatibility
between $*$-structures and
involutions with the coalgebra
structure:
\begin{definition}
Let $H$ be 
a Hopf algebra. 
\begin{enumerate}
\item A \emph{Hopf
$*$-structure} on $H$ is a
$*$-structure on the
underlying algebra satisfying 
for all $h \in H$ 
$$ 
	(h^*)_{(1)} \otimes (h^*)_{(2)} 
	= 
	 (h_{(1)})^* \otimes (h_{(2)})^*,
	\quad \varepsilon (h^*) = \overline{
	\varepsilon (h)}.
$$ 
\item A \emph{Cartan
involution} on
$H$ is an involution 
$\theta$ on the underlying
algebra such that for all 
$h \in H$, we have
$$
	\theta(h) _{(1)} \otimes 
	\theta(h) _{(2)} =  
	\theta(h_{(2)}) \otimes 
	\theta(h_{(1)}),\quad
	\varepsilon (\theta(h)) 
	=\overline{\varepsilon
(h)}.
$$
\end{enumerate}
\end{definition}

These structures correspond
bijectively to each other:

\begin{lemma}\label{starcartan}
A map $* \colon H \rightarrow
H$ is a Hopf
$*$-structure if and only if 
$* \circ S$ is a Cartan
involution.
\end{lemma}
\begin{proof}
This is mostly obvious; for a
proof of 
$ * \circ S \circ *
\circ S= \mathrm{id} _H$ 
that carries over verbatim to
arbitrary ground fields, see
e.g.~\cite[Proposition~1.2.7.10]{schmudgen}.
\end{proof}
\begin{remark}
More abstractly, a Cartan
involution corresponds to
a lift of the endofunctor $V \mapsto
\overline{V}$ from the category of 
finite-dimensional $k$-vector
spaces to that of
finite-dimensional 
$H$-modules: since $H$
is a Hopf algebra, this
category is a
rigid monoidal category, that is, one can
form tensor products and
left and right duals of $H$-modules, and
these structures are compatible with the
forgetful functor to $k$-vector spaces. 
A Hopf
$*$-structure on $H$ 
allows one to additionally define for any
$H$-module $V$ 
the conjugate
$\overline{V}$ to be the conjugate vector space
with $H$-module structure given by 
$h v:=\thetah(h) v=S(h)^*v$.
Note that as $\theta$ is a
coalgebra morphism 
$ H \rightarrow \bar H^\cop$,
the compatibility with
the monoidal structure 
is $\overline{V \otimes W} = 
\overline{W} \otimes \overline{V}$. 
\end{remark}

If $H$
is a Hopf $*$-algebra, then 
an $H$-module algebra $B$ 
which is also a $*$-algebra is
called a \emph{module
$*$-algebra} if 
$(hb)^* = {\thetah}(h)(b^*)$ 
holds for all $h \in
H,b \in B$, that is, if the
resulting map $H \rightarrow
\mathrm{End} _k(B)$ is a
morphism of algebras with
involution. The question we
want to address in this
subsection 
is whether for a given
quantum automorphism $ \sigma
$ of a $*$-algebra $B$ the
Hopf algebra $H_\sigma$
becomes naturally a
Hopf $*$-algebra in such a way that
$B$ is a module $*$-algebra
over $H_\sigma$. 
 
In order to do so, we first 
extend the $*$-structure $*$ respectively
the associated involution ${\thetab}$ 
to $\Mn(B)$ respectively 
$\Mn( \mathrm{End}_k(B))$. This depends
on the choice of an involutive permutation 
$$
	\{1,\ldots,n\} \rightarrow 
	\{1,\ldots,n\},\quad
	i \mapsto \bar i,\quad
	\bar{\bar{i}}=i
$$
that will be used afterwards to
be able to restrict the resulting
$*$-structure on $\Mn(B)$ to upper
triangular matrices. 

\begin{proposition}\label{starstr}
Let $B$ be a $*$-algebra and 
assume $s \in S_n$, 
$s^2=1$. We abbreviate $\bar i:=s(i)$. 
\begin{enumerate}
\item Setting 
$
	(\sigma^\dagger)_{ij} := 	
	\sigma_{\bar j \bar i}^* 
$
defines a $*$-structure $\dagger$ on $\Mn(B)$.
\item Setting 
$
	\vartheta (\sigma)_{ij} := 	
	{\thetab}(\sigma_{\bar i \bar j})  
$
defines an involution on $\Mn(
\mathrm{End}_k(B))$.  
\end{enumerate}
\end{proposition}
\begin{proof}
Clearly, $\dagger$ is involutive:
$$
	(\sigma^{\dagger\dagger})_{ij} =
(\sigma^\dagger)^*_{\bar j \bar i} 
	= \sigma^{**}_{\bar{\bar i} \bar{\bar j}} 
	= \sigma^{**}_{ij} = \sigma_{ij},
$$
where the last equality
follows as 
$*$ is a $*$-structure on $B$. It is also a
ring morphism $\Mn(B) \rightarrow
\Mn(B)^\op$: let
$\sigma, \tau \in \Mn(B)$ then
\begin{align*}
		((\sigma \tau)^\dagger)_{ij} 
&= ((\sigma \tau)_{\bar j \bar i})^* 
			= \big(\sum_{r = 1}^n  \sigma_{\bar j r}\tau_{r \bar i}\big)^* 
			= \sum_{r = 1}^n  \tau^*_{r \bar i}\sigma^*_{\bar j r} \\
			&= \sum_{r = 1}^n
(\tau^\dagger)_{i \bar r }(\sigma^\dagger)_{\bar r j} 
			= (\tau^\dagger \sigma^\dagger)_{ij}.
\end{align*}
That $\dagger$ is a $k$-linear map $\Mn(B)
\rightarrow \overline{\Mn(B)}$ follows from
the fact that $* \colon B \rightarrow \bar
B$ is linear.
The second claim is shown analogously.  
\end{proof}

The definition of this $*$-structure and
of this involution is made in order to
have the following: 

\begin{lemma}
A $k$-linear map $ \sigma \colon B \rightarrow \Mn(B)$ 
satisfies 
$$
	\sigma (b^*) = \sigma (b)^\dagger
$$
if and only if $ \vartheta( \sigma ) =
\sigma^T $ when $ \sigma $ is
viewed as an element in $
\Mn(\mathrm{End}_k(B))$.  
\end{lemma}
\begin{proof}
This holds as
$(\sigma^\dagger)_{ij}(b) 
=\sigma^*_{\bar j \bar i}(b)
= {\thetab}(\sigma_{\bar j \bar i})(b^*)
=\vartheta(\sigma)_{ji}(b^*)$.
\end{proof}

Now we apply the above to the study of
quantum automorphisms:

\begin{proposition}
Let $ \sigma \colon B \rightarrow \Mn(B)$ be a
quantum automorphism. 
\begin{enumerate}
\item 
$ \vartheta(\sigma)^T=\vartheta(\sigma ^T)$ is a quantum
automorphism.
\item If $\vartheta( \sigma )^T = \sigma
$, then $H_\sigma $ is a Hopf $*$-algebra
with $*$-structure given by
$
	s_{d,ij} \mapsto 
	s_{1-d,\bar i\bar j},
$
and $B$
is a module $*$-algebra over $H_\sigma$. 
\end{enumerate}
\end{proposition}
\begin{proof}
(1): Since $ \sigma \colon B \rightarrow \Mn(B)$
is an algebra morphism, we have
\begin{align*}
	(\vartheta (\sigma)^T)_{ij} (ab)  
&=
	{\thetab} (\sigma_{\bar j \bar i}) (ab)  
= \sigma_{\bar j \bar i} ((ab)^*) ^*
= \sigma_{\bar j \bar i} (b^*a^*)^* \\
&= (\sum_r \sigma _{\bar j r}
	(b^*) \sigma_{r \bar i}
	(a^*))^* \\
&= \sum_r 
	(\sigma_{r \bar i}
	(a^*))^* 
	(\sigma _{\bar j r}
	(b^*))^* 
\\
&= \sum_r 
	{\thetab}(\sigma_{r \bar i}) 
	(a) 
	{\thetab}(\sigma _{\bar j r})
	(b) 
\\
&= \sum_r
	(\vartheta(\sigma)^T)_{i\bar r}(a) 
	(\vartheta (\sigma)^T)_{\bar rj}(b) \\
&= \sum_r
	(\vartheta(\sigma)^T)_{ir}(a) 
	(\vartheta (\sigma)^T)_{rj}(b).
\end{align*}	
In order to 
show that 
$ \vartheta( \sigma )^T$ is strongly
invertible, note first that 
as $ \vartheta $ is an 
involution on $\Mn(
\mathrm{End}_k(B) )$,  it is in
particular multiplicative, 
$$
	\vartheta ( \sigma \tau ) = 
	\vartheta (\sigma ) 
	\vartheta ( \tau ),
$$ 
so if $ \sigma \in 
\Mn( \mathrm{End}_k(B)) $ is
invertible, then so is 
$ \vartheta ( \sigma )$ with 
inverse given by 
$ \vartheta ( \sigma )^{-1} = 
\vartheta( \sigma ^{-1})$. 
Furthermore, it follows
directly from its definition
that $\vartheta$ commutes with
taking transposes, so 
$$
	\vartheta (\hat \sigma ) 
	=\vartheta ((\sigma
^{-1})^T) = 
	\vartheta (\sigma ^{-1}) ^T
= 
	(\vartheta (\sigma) ^{-1}) ^T
= 
	\widehat {\vartheta
(\sigma )}.
$$	
It follows that if 
$ \{\sigma _d\}$ is a sequence
of invertible matrices with 
$ \hat \sigma _d = \sigma
_{d+1}$ then 
$ \{\vartheta (\sigma_d) \}$
is a sequence of invertible
matrices with 
$ \widehat{\vartheta(\sigma
_d) } = \vartheta (\sigma
_{d+1})$. 
(2): Since $k \langle s_{d,ij}
\rangle $ is a free algebra, 
there is a unique
algebra morphism 
$$
	\theta \colon 
	k \langle s_{d,ij} \rangle 
	\rightarrow \overline 
	{k \langle s_{d,ij}
\rangle},\quad
	s_{d,ij} \mapsto  
	s_{-d,\bar j\bar i}.
$$ 
That $ k \langle s_{d,ij} \rangle 
\rightarrow \overline{
k \langle s_{d,ij} \rangle
}^\cop$ is a coalgebra morphism 
 is verified by straightforward
computation.

Our aim is to show that $ \theta $ descends to a Cartan 
involution on $H_\sigma$. 
For this we first prove by induction on
$d$ that 
\begin{equation}\label{drill}
	\thetab( \sigma _{d,\bar i \bar j}) = 
	\vartheta ( \sigma _d)_{ij} = 
	\sigma _{-d,ij}^T = 
	\sigma _{-d,ji}.
\end{equation}
Indeed, for $d=0$ this holds by
the assumption that $ \vartheta ( \sigma
)^T=\sigma $. In the induction step we
compute  
$$
	\vartheta (\sigma_d ) = \sigma_{-d}^T 
	\Rightarrow 
	\vartheta (\sigma_d ^{-1}) =
(\sigma_{-d}^T)^{-1} = \sigma _{-d-1}, 
$$
where the last equality
follows from the comment
after 
Definition~\ref{StrInv}.
Hence
$$
	\vartheta ( \sigma _{d+1}) =
	\widehat{\vartheta(\sigma_d )} = 
	\vartheta (\hat \sigma_d ) = 
	\vartheta ( \sigma _d^{-1}) ^T = 
	\sigma _{-d-1}^T
$$
which proves (\ref{drill}) for $d+1$. 
This means that by the definition of the
involution $\theta$ on $k \langle s_{d,ij}
\rangle $, the map $ k \langle s_{d,ij}
\rangle \rightarrow \mathrm{End}_k(B)$,
$s_{d,ij} \mapsto \sigma _{d,ij}$ is a
morphism of algebras with involution. In
particular, $ \theta $ descends to 
$k \langle s_{d,ij} \rangle /I$ and to 
$ H_ \sigma $, and defines a Cartan
involution and hence by
Lemma~\ref{starcartan} a Hopf $*$-algebra
structure on $H_\sigma$. It also follows
that $B$ is a module $*$-algebra over this
Hopf $*$-algebra. 
\end{proof}

\begin{example}\label{invperms}
Throughout, we fixed an involutive
permutation $s \in S_n$. The most obvious
choice is the identity, $ \bar i=i$. In
this case, $ \sigma ^\dagger $ is the
usual adjoint of a matrix $ \sigma \in
\Mn(B)$ (transpose and apply $*$ to the
entries). However, if $ \sigma $ is upper
triangular, the condition 
$ \vartheta ( \sigma )^T= \sigma $ can not
hold with respect to this involution, as
we then have $\vartheta ( \sigma )_{ij} = 
\thetab(\sigma _{ij}) $, so
that $\vartheta
(\sigma )^T$ is a lower triangular matrix. 
Hence for upper
triangular quantum
automorphisms we focus 
on the permutation 
$\bar i := n+1-i$.  

In particular, consider $k= \mathbb{C} $ 
with involution given by
complex conjugation and
$B=\mathbb{C} [t^2,t^3]$ with
$*$-algebra structure given by  
$$
	* \colon B= \mathbb{C} [t^2,t^3] \longrightarrow B, 
	\quad (\lambda t^n)^* := 
	\bar \lambda t^n.
$$
Geometrically, this 
describes the real points of the singular
curve $V \subset 
\mathbb{C}^2$ in the sense that the points of
the curve $V_\mathbb{R} = V \cap
\mathbb{R} ^2$ correspond to the
one-dimensional $*$-representations of
$B$, and these correspond
further to the maximal
ideals in $ \mathbb{C} [t^2,t^3] $ which
are invariant under $*$. 

The quantum automorphism 
$ \sigma \colon B \rightarrow \Mthree(B)$
that we study since Section~\ref{explicit} 
satisfies $\vartheta ( \sigma )^T =
\sigma $ provided that we work with the
involution  
$\bar j := 4-j$ on $\{1,2,3\}$ as in
Example~\ref{invperms}. The Hopf
$*$-structure on $H_\sigma$ is then given
by 
$$
	K^* = K,\quad
	D^* = -D,\quad
	Y^* = -Y + 6iD. 
$$
\end{example}


\end{document}